\documentclass[12pt,twoside,a4paper]{article}
\usepackage{stmaryrd}
\usepackage{amsfonts}
\usepackage{amsmath,amssymb}
\usepackage{mathrsfs}
\usepackage{hyperref}

\newtheorem{thm}{Theorem}[section]
\newtheorem{cor}[thm]{Corollary}
\newtheorem{lem}[thm]{Lemma}
\newtheorem{prop}[thm]{Proposition}

\newtheorem{rem}[thm]{Remark}

\numberwithin{equation}{section}\allowdisplaybreaks

\def\le{\leqslant}
\def\ge{\geqslant}
\def\leq{\leqslant}
\def\geq{\geqslant}
\def\no{\nonumber}

\begin{document}

\begin{center}
\Large\bf  Analyticity for the (generalized) Navier-Stokes equations with rough initial data

 \vspace*{0.5cm}\normalsize
  Chunyan Huang$^{a}$,\quad Baoxiang Wang$^{b,}$\footnote{Corresponding author}

\vspace*{0.5cm} {\it \footnotesize $^{a}$School of Statistics and
Mathematics, Central University of Finance
and Economics, Beijing 100081 People's Republic of China, Email:
  hcy@cufe.edu.cn\\ $^{b}$LMAM, School of Mathematical Sciences, Peking University,
Beijing 100871 People's Republic of China, Email:
  wbx@math.pku.edu.cn}\\

\vspace*{0.5cm}

\begin{minipage}{13.5cm}
\footnotesize \bf Abstract. \rm  We study the Cauchy problem for the (generalized)  incompressible
Navier-Stokes equations
\begin{align}
u_t+(-\Delta)^{\alpha}u+u\cdot \nabla u +\nabla p=0, \ \  {\rm div} u=0, \ \  u(0,x)= u_0. \nonumber
\end{align}
We show the analyticity of the local solutions of the Navier-Stokes equation ($\alpha=1$) with any initial data  in critical Besov spaces $\dot{B}^{n/p-1}_{p,q}(\mathbb{R}^n)$ with  $1< p<\infty, \ 1\le q\le \infty $  and the solution is global if $u_0$ is sufficiently small in $\dot{B}^{n/p-1}_{p,q}(\mathbb{R}^n)$. In the case $p=\infty$, the analyticity for the local solutions of the Navier-Stokes equation ($\alpha=1$) with any initial data in modulation space $M^{-1}_{\infty,1}(\mathbb{R}^n)$ is obtained.  We prove the
global well-posedness for a fractional Navier-stokes
equation ($\alpha=1/2$) with small data in critical Besov spaces
$\dot{B}^{n/p}_{p,1}(\mathbb{R}^n) \ (1\leq p\leq\infty)$ and show the
analyticity of solutions with small initial data either in $\dot{B}^{n/p}_{p,1}(\mathbb{R}^n) \ (1\leq p<\infty)$ or in $\dot{B}^0_{\infty,1} (\mathbb{R}^n)\cap  {M}^0_{\infty,1}(\mathbb{R}^n)$. Similar results also hold for all $\alpha\in (1/2,1)$. \\

\\

\end{minipage}
\end{center}

\section{Introduction} \label{sect1}

\medskip

In this paper, we study the Cauchy problem for the (generalized)  incompressible
Navier-Stokes equations (GNS):
\begin{align}
u_t+(-\Delta)^{\alpha}u+u\cdot \nabla u +\nabla p=0, \ \  {\rm div} u=0, \ \  u(0,x)= u_0. \label{1.1}
\end{align}
where $t\in \mathbb{R}^{+}=[0,\infty)$, $x\in \mathbb{R}^{n} \ (n\geq 2)$, $\alpha
> 0$ is a real number. $u=(u_1,...,u_n)$ denotes the flow
velocity vector and $p(t,x)$ describes the scalar pressure.
$u_0(x)=(u_0^1,...,u_0^n)$ is a given velocity with $\mathrm{div}
u_0=0$. $u_t=\partial u/\partial t,$ $(-\Delta)^{\alpha}$ denotes
the fractional Laplacian which is defined as
$\widehat{(-\Delta)^{\alpha}}u(t,\xi)=|\xi|^{2\alpha}\widehat{u}(t,\xi)$.
It is easy to see that \eqref{1.1} can be rewritten as the following
equivalent form:
\begin{align}
  u_t+(-\Delta)^{\alpha}u+\mathbb{P}\ \textrm{div}(u\otimes u)=0, \ \
u(0, x)= u_0,
 \label{2.1}
\end{align}
where $\mathbb{P}=I-{\rm div} \Delta^{-1}\nabla$ is the matrix operator
projecting onto the divergence free vector fields,  $I$ is identity
matrix.

It is known that \eqref{1.1} is essentially equivalent to the following integral equation:
\begin{align} \label{1.3}
u(t) = U_{2\alpha}(t) u_0 +  \mathscr{A}_{2\alpha}\mathbb{P}\ \textrm{div}(u\otimes u),
\end{align}
where
\begin{align}
U_{2\alpha} (t):=e^{-t(-\Delta)^{\alpha}}=\mathscr{F}^{-1}e^{-t|\xi|^{2\alpha}}\mathscr{F}, \quad   (\mathscr{A}_{2\alpha}f)(t):= \int_0^t U_{2\alpha}(t-\tau)f(\tau)d\tau.
\end{align}

Note that \eqref{1.1} is scaling invariant in the following sense:
if $u$ solves \eqref{1.1},  so does   $u_{\lambda}(t, x)=\lambda^{2\alpha-1}u(\lambda^{2\alpha}t, \lambda x)$ with initial data $\lambda^{2\alpha-1}u_0(\lambda x)$.
 A function space $X$ defined in $\mathbb{R}^n$ is said to be a critical space for \eqref{1.1}
if the corresponding norm of $u_{\lambda}(0,x)$ in $X$ is invariant for
all $\lambda>0$. It is easy to see that homogeneous type Besov space
$\dot{B}^{n/p-2\alpha +1}_{p, q}$ is a critical space for
\eqref{1.1}.  In particular, $\dot{B}^{n/p-1}_{p, q}$ is a critical space for GNS in the case $\alpha=1$ and $\dot{B}^{n/p}_{p, q}$ is a critical space for GNS in the case
$\alpha=1/2$.

If $\alpha=1$, GNS \eqref{1.1} is the well-known
Navier-Stokes equation (NS) which have been extensively studied  in the past twenty years; cf. \cite{BaBiTa12, Can97, CaPl96, Ch99, DoLi09, EsSeSv03, FoTe89, GePaSt07, GrKu98, Iw10, Ka84, KoTa01, LeRi00, Pl96}. Kato \cite{Ka84} first used the semi-group estimates to obtain the strong $L^p$ solutions of NS. Cannone \cite{Can97} and Planchon \cite{Pl96} considered global solutions in 3D for small data  in critical Besov spaces $\dot{B}^{3/p-1}_{p, \infty}$  with $3<p\le 6$. Chemin \cite{Ch99} obtained global solutions in 3D for small data  in critical Besov spaces $\dot{B}^{3/p-1}_{p, q}$ for all $p<\infty, \ q\le \infty$.  Koch and Tataru  \cite{KoTa01} studied local solutions for initial data in $\mathit{vmo}^{-1}$ and global solutions for small initial data in $BMO^{-1}$. Escauriaza, Serigin and  Sverak were able to show the 3D global regularity if the solution has a uniform bound in the critical Lebesgue space $L^3$. Iwabuchi \cite{Iw10} considered the well posedness of NS for initial data in modulation spaces $M^{s}_{p,q}$, especially in $M^{-1}_{p,1}$ with $1\le p\le \infty$. Foias and Temam \cite{FoTe89} proved spatial analyticity for solutions in Sobolev spaces of periodical functions in an elementary way. The analyticity of solutions in $L^p$ for NS was first shown by Gruji\v{c} and Kukavica \cite{GrKu98} and Lemari\'{e}-Rieusset \cite{LeRi00} gave a different approach based on multilinear singular integrals. Using iterative derivative estimates, the analyticity of NS for small initial data in $BMO^{-1}$ was obtained in Germain, Pavlovic and Staffilani \cite{GePaSt07} (see also Dong and Li \cite{DoLi09}, Miura and Sawada \cite{MiSa06} on the iterative derivative techniques).
Recently, Bae, Biswas and  Tadmor \cite{BaBiTa12} obtained the analyticity of the solutions of NS in 3D for the sufficiently small initial data in critical Besov spaces $\dot{B}^{3/p-1}_{p, q}$ with $1< p,q<\infty$. Moreover, in the case $p=\infty$, they can show the analyticity of the solutions of NS in 3D for the sufficiently small initial data in $\dot{B}^{-1}_{\infty, q} \cap \dot{B}^{0}_{3, \infty}$ with $1\le q<\infty$. Since NS is ill-posed in $\dot{B}^{-1}_{\infty,\infty}$; cf. Bourgain and Pavlovic \cite{BoPa08}, Germain \cite{Ge08} and Yoneda \cite{Yo10}, the largest Besov-type space on initial data for which NS is well-posed or analytically well-posed is still open.

 For general $\alpha$,  there are
also some results on \eqref{1.1} in recent years. For $\alpha \geq
5/4,  n=3$,  Lions \cite{Li} proved the global existence of classical
solutions.  Wu \cite{WU3} proved \eqref{1.1} has a unique global
solution with small data in $\dot{B}^{n/p-2\alpha+1}_{p, q} \ (1\leq q
\leq\infty)$ for $\alpha
>1/2$,  $p=2$,  or $1/2<\alpha\leq1$,  $2<p<\infty$. Recently,  Li, Zhai and
Yang \cite{LZ, LY} obtained the well-posedness of \eqref{1.1} in the case
$1/2<\alpha<1$ in some Q-spaces. Yu and Zhai \cite{YuZh12} showed the global existence and uniqueness with small initial data in $\dot{B}^{1-2\alpha}_{\infty,\infty}$ for GNS in the case $1/2<\alpha <1$.  As far as the authors can see,  for $\alpha =1/2$ and general $p$,  the
well-posedness and analyticity of solutions of \eqref{1.1} in critical Besov spaces
$\dot{B}^{n/p}_{p, q}$ seems unsolved.

\subsection{Notations}\label{sect2}

Throughout this paper, $C\ge 1, \ c\le 1$ will denote constants which can be different at different places, we will use $A\lesssim B$ to denote   $A\leqslant CB$. We denote  $|x|_p=(|x_1|^p+...+|x_n|^p)^{1/p}$, $|x|=|x|_2$ for any $x\in \mathbb{R}^n$,  by $L^p=L^p(\mathbb{R}^n)$ the Lebesgue space on which the norm is written as $\|\cdot\|_p$.  Let $X$ be a quasi-Banach space. For any $I\subset \mathbb{R}^+=[0,\infty)$,   we denote
$$
\|u\|_{ L^{\gamma}(I; \ X)} = \left(\int_I \|u(t, \cdot)\|^{\gamma}_{X} dt \right)^{1/{\gamma}}
$$
 for $1\leq \gamma <\infty$ and with usual modifications
for $\gamma =\infty$. In particular,  if $X=L^p$,  we will write
$\|u\|_{L_{t\in I}^{\gamma}L_x^p}=\|u\|_{L^{\gamma}(I;L^p)}$ and $\|u\|_{L_{t }^{\gamma}L_x^p}=\|u\|_{L^{\gamma}(0,\infty;L^p)}$. For any function space $X$ and the operator $T: X\to X,$  we denote
\begin{align}
T X=\{T f : \ f\in X\}, \quad \|f\|_{T X} := \|Tf\|_X.  \label{1.5}
\end{align}
Now let us recall the definition of dyadic decomposition in
Littlewood-Paley theory \cite{Tr}. Let $\psi: \mathbb{R}^n\rightarrow
[0, 1]$ be a smooth cut-off function which equals $1$ on the unit ball and equals $0$ outside the ball $\{\xi
\in \mathbb{R}^n; |\xi|\leq 2\}$.
Write $\varphi(\xi):=\psi(\xi)-\psi(2\xi)$ and
$\varphi_k(\xi)=\varphi(2^{-k}\xi)$.
$\Delta_k:=\mathscr{F}^{-1}\varphi_k \mathscr{F}, k\in \mathbb{Z}$
are said to be the dyadic decomposition operators and satisfying the
operator identity: $I=\sum_{k=-\infty}^{+\infty}\Delta_k$.
$S_k:=\sum_{j=-\infty}^{k}\Delta_j$ are called the low frequency
projection operators. It is easy to see that $S_k u \rightarrow u$
as $k\rightarrow \infty$ in the sense of distributions. The norms in homogeneous Besov spaces are defined as follows:
\begin{align}
\|f\|_{\dot{B}^s_{p, q}}=
 \left(\sum_{j=-\infty}^{+\infty}2^{jsq}\|\Delta_j f\|^q_{p}\right)^{1/q}
\label{Besov}
\end{align}
with usual modification if $q=\infty.$ We also need the following  time-space norms:
\begin{align} \label{LBesov}
\|f\|_{\widetilde{L}^{\gamma}(I;\dot{B}^s_{p, q})} =
\left(\sum_{j=-\infty}^{+\infty}2^{jsq}\|\Delta_j
f\|^q_{L_{t\in I}^{\gamma}L_x^p} \right)^{1/q}<\infty
\end{align}
with the usual modification for $q=\infty$ and by \eqref{1.5}, for any $ T :\dot{B}^s_{p, q} \to \dot{B}^s_{p, q}$,
\begin{align}
\|f\|_{\widetilde{L}^{\gamma}(I; T \dot{B}^s_{p, q})}
= \|T f\|_{\widetilde{L}^{\gamma}(I;\dot{B}^s_{p, q})},
\end{align}
Now
we recall the definition of modulation spaces which was first
introduced by Feichtinger \cite{Fei83} in 1983 (see also Gr\"ochenig \cite{Groch01}). Let $\sigma$ be a smooth
cut-off function with ${\rm supp} \ \sigma \subset  [-3/4, 3/4]^n$,
$\sigma_k=\sigma(\cdot-k)$ and
 $  \sum_{ k\in\mathbb{Z}^n}\sigma_k = 1.$
We define the frequency-uniform decomposition operator as
   $\square_k = \mathscr{F}^{-1} \sigma_k \mathscr{F},  \quad \sigma \in \mathbb{Z}^n. $
Let $0<p, q\le \infty$,  $s\in \mathbb{R}$ and
  \begin{align}
\|f\|_{M^s_{p, q}} =  \left(\sum_{k\in \mathbb{Z}^n}  \langle k\rangle^{sq} \|\Box_k   f\|^q_{p} \right)^{1/q},   \label{modsp}
\end{align}
$M^s_{p, q}$ is said to be a modulation space. There are exact inclusions between Besov spaces and modulation spaces (see \cite{Gr92, To04, SuTo07, WaHu07}), however, there are no straightforward inclusions between homogeneous Besov spaces and modulation spaces.  Let $0<p, q\le \infty$,  $s\ge 0$,
\begin{align}
\|f\|_{E^s_{p, q}} =  \left(\sum_{k} 2^{q s | k |} \|\Box_k   f\|^q_{p} \right)^{1/q},   \label{emodsp}
\end{align}
$E^s_{p, q}$ ($s>0$) was introduced in \cite{Wa06}, which can be regarded as modulation spaces with analytic regularity. It is easy to see that $M^0_{p,q}= E^{0}_{p,q}$.  Let $\alpha=(\alpha_1, ..., \alpha_n)$,  $\alpha! = \alpha_1!...\alpha_n!$  and $\partial^\alpha=\partial^{\alpha_1}_{x_1}... \partial^{\alpha_n}_{x_n}$.
Recall that the Gevrey class is defined as follows.
$$  G_{1, p}= \left\{f\in C^\infty(\mathbb{R}^n):\,  ^\exists\rho,  M>0\; s.t.\;
 \|\partial^\alpha f \|_p \le   \frac{M \alpha!}{\rho^{|\alpha|}}  ,  \ \forall \
 \alpha\in {\Bbb
Z}^n_+  \right\}.
$$
It is known that $G_{1,\infty}$ is the Gevrey 1-class. Moreover, we easily see that $G_{1,p_1} \subset G_{1, p_2}$ for any $p_1\le p_2$.  There is a very beautiful relation between Gevrey class and exponential modulation spaces, we can show that  $G_{1,p}=\bigcup_{s>0}E^s_{p, q}$, see Section \ref{sectGevExp}.
Similarly as in \eqref{LBesov}, we need the following norm
\begin{align}
\|f\|_{\widetilde{L}^{\gamma}(I; {E}^{s(t)}_{p, q})}=
\left(\sum_{k\in \mathbb{Z}^n}\|2^{ s(t) |k|} \Box_k
f\|^q_{L_{t\in I}^{\gamma}L_x^p} \right)^{1/q},
\end{align}
with usual modification if $q=\infty$. Recently modulation spaces have been applied to the study of nonlinear evolutions; cf. \cite{BenOk08,BGOR, ChFaSu12,CorNik08,CorNik092,Iw10,WaHuHaGu11,Wa06}.

As the end of this section, we recall the definition of multiplier space $M_p$ (see \cite{BL, Ho60}).
Let $\rho \in \mathscr{S}'$. If there exists a $C>0$ such that
$\|\mathscr{F}^{-1}\rho \mathscr{F}f\|_{L^p}\leq C\|f\|_{L^p}$ holds
for all $f\in \mathscr{S}$, then $\rho$ is called a Fourier
multiplier on $L^p$. The linear space of all multipliers on $L^p$ is
denoted by $M_p$ and the norm on which is defined as
$\|\rho\|_{M_p}=\sup\{\|\mathscr{F}^{-1}\rho \mathscr{F}f\|_{L^p}:
f\in \mathscr{S}, \|f\|_{L^p}=1\}$. Using the multiplier estimates, we have the following Bernstein's inequalities:
\begin{align*}
\|\Delta_j (-\Delta)^{s/2} f\|_{q} & \lesssim
2^{j(s+n(1/p-1/q))}\|f\|_{p}, \\
\|\Delta_j (-\Delta)^{s/2} f\|_{p} & \sim 2^{js}\|\Delta_j f\|_{p}
\end{align*}

\subsection{Main Results}

The analyticity for the global solutions of NS in 3D with small initial data in critical Besov spaces $\dot{B}^{3/p-1}_{p,q}$ ($1\le p,q<\infty$) was obtained in \cite{BaBiTa12}. In this paper, we show the analyticity of the local solutions of NS with large initial data in critical Besov spaces $\dot{B}^{n/p-1}_{p,q}$ ($1< p<\infty, \ 1\le q\le \infty$). Let us denote $\Lambda=(-\partial^2_{x_1})^{1/2}+...+(-\partial^2_{x_n})^{1/2}$. Recall that $\|f\|_{\widetilde{L}^{q}(I;  e^{\sqrt{t}\Lambda}\dot{B}^{s}_{p, q})} = \|e^{\sqrt{t}\Lambda} f\|_{\widetilde{L}^{q}(I;\dot{B}^{s}_{p, q})}$. Our main results are the following theorems.

\begin{thm}\label{NSTh1a} {\rm (Analyticity for NS (I): $p<\infty$)}
Let $n\ge 2$, $ \alpha=1$, $1< p < \infty$, $1\le q\le \infty$.  Assume that $u_0 \in   \dot{B}^{n/p-1}_{p, q}$.  There exists a $T_{\max}= T_{\max}(u_0)>0$ such that \eqref{1.1} has a unique solution
$$
u \in  \widetilde{L}^{3}_{\rm loc}(0,T_{\max};  e^{\sqrt{t}\Lambda}\dot{B}^{n/p-1/3}_{p, q}) \cap \widetilde{L}^{3/2}_{\rm loc}(0,T_{\max};  e^{\sqrt{t}\Lambda}\dot{B}^{n/p+1/3}_{p, q}) \cap  \widetilde{L}^{\infty}([0,T_{\max}); e^{\sqrt{t}\Lambda} \dot{B}^{n/p-1}_{p, q} ).
$$
If $T_{\max}<\infty$, then
$$
  \|u\|_{\widetilde{L}^{3}(0,T_{\max};  e^{\sqrt{t}\Lambda}\dot{B}^{n/p-1/3}_{p, q}) \cap \widetilde{L}^{3/2}(0,T_{\max};  e^{\sqrt{t}\Lambda}\dot{B}^{n/p+1/3}_{p, q})}=\infty.
$$
 Moreover, if $u_0 \in   \dot{B}^{n/p-1}_{p, q}$ is sufficiently small, then $T_{\max}=\infty.$
\end{thm}

Theorem \ref{NSTh1a} has implied the analyticity of the solutions. In the case $p=\infty$, the analyticity of global solutions of NS with small initial data in $\dot{B}^{-1}_{\infty, q} \cap \dot{B}^{0}_{3, \infty}$ was shown in \cite{BaBiTa12}.  We use a different approach to get an analyticity result for local solutions with any initial data in the modulation space $M^{-1}_{\infty,1}$:

\begin{thm}\label{NSTh1b} {\rm (Analyticity for NS (II): $p\le \infty$)}
Let $n\ge 2$, $ \alpha=1$, $1\le p \le \infty$.  Assume that $u_0 \in   M^{-1}_{p, 1}$.  There exists a $T_{\max}= T_{\max}(u_0)>0$ such that \eqref{1.1} has a unique solution $u \in  \widetilde{L}^{2}_{\rm loc}(0,T_{\max}; {E}^{ct}_{p, 1})$ and $ (I-\Delta)^{-1/2} u \in  \widetilde{L}^{\infty}([0,T_{\max}); {E}^{ct}_{p, 1})$, $ \ c=2^{-10} $. Moreover,
 if $T_{\max}<\infty$, then $  \|u\|_{\widetilde{L}^{2}(0,T_{\max}; {E}^{ct}_{p, 1})}=\infty$.
\end{thm}

Now we compare homogeneous Besov spaces with modulation spaces. Let $S_1$ be the low frequency projection operator defined in the above. Let $1\le p\le \infty$.
One easily sees that $S_1 \dot{B}^{-1}_{p,1} \subset S_1 M^{-1}_{p,1}$ and $(I-S_1) M^{-1}_{p,1} \subset (I-S_1) \dot{B}^{-1}_{p,1}$. Roughly speaking, the low frequency part of  $M^{-1}_{p,1}$ is rougher than that of $\dot{B}^{-1}_{p,1}$, the high frequency part of $M^{-1}_{p,1}$ is smoother than that of $\dot{B}^{-1}_{p,1}$. In general there is no inclusion between $M^{-1}_{p,1}$ and $\dot{B}^{-1}_{p,1}$.

If $\alpha=1/2,$ we have

\begin{thm}\label{Th1}
Let $n\ge 2$,  $ \alpha=1/2$.  Suppose $u_0 \in \dot{B}^{n/p}_{p, 1} \ (1\leq
p  \le \infty)$ and $\|u_0\|_{\dot{B}^{n/p}_{p, 1}}$  is sufficiently
small,  then we have the following results:
\begin{itemize}
\item[\rm i)] {\rm (Well-posedness for $p<\infty$)} If $1\le p <\infty$, then there exists a unique global solution $u$ to
\eqref{1.1} satisfying
$u\in
\widetilde{L}^{\infty}(\mathbb{R}_{+};\dot{B}^{n/p}_{p, 1})$.

\item[\rm ii)] {\rm (Well-posedness for $p=\infty$)} If $p=\infty$, then there exists a unique global solution $u$ to
\eqref{1.1} satisfying $u\in
\widetilde{L}^{\infty}(\mathbb{R}_{+};\dot{B}^{0}_{\infty,1})\cap \widetilde{L}^1(\mathbb{R}_{+},\dot{B}^{1}_{\infty,1})$.

\item[\rm iii)] {\rm (Analyticity  for $p< \infty$)} If $1< p <\infty$, then the solution obtained in {\rm i)} satisfying  $u(t) \in \widetilde{L}^{\infty}(\mathbb{R}^+,  e^{t \Lambda/2n} \dot{B}^{n/p}_{p,1})$.

\item[\rm iv)] {\rm (Analyticity  for $p=\infty$)}
  Assume that $u_0 \in \dot{B}^{0}_{\infty,1} \cap M^0_{\infty, 1}$ and $\|u_0\|_{\dot{B}^{0}_{\infty,1} \cap M^0_{\infty, 1}}$ is sufficiently small.   Then the solution obtained in {\rm ii)} satisfies  $u \in \widetilde{L}^{\infty}(\mathbb{R}^+,  E^{s(t) }_{\infty, 1})$, $\nabla u \in \widetilde{L}^{1}(\mathbb{R}^+,  E^{s(t) }_{\infty, 1})$ for $s(t)=2^{-5} (1\wedge t)$.

\end{itemize}
\end{thm}
In Theorem \ref{Th1}, the uniqueness of the well-posedness of the case $p=\infty$ is worse than that of the case $p<\infty$.  Moreover, the proof of Theorem \ref{Th1} cannot be developed to the case for which the initial data belong to the other critical Besov spaces $\dot{B}^{n/p}_{p,q}$ with $q\neq 1$.

At the end of this section, we point out that the ideas used in this paper is also adapted to the case $1/2<\alpha <1$. In Section \ref{GNS7} we will generalize the analyticity results  in Theorems \ref{NSTh1a} and \ref{NSTh1b} to the case $1/2<\alpha <1$, see Theorems \ref{GNSTh1a} and \ref{GNSTh1b}.

\section{Analyticity of NS: Proof of Theorem \ref{NSTh1a}} \label{sectNSa}

First, we establish two linear estimates by following some ideas as in \cite{Ch99, LeRi02, BaBiTa12} and \cite{Wang04,Wang02}.
If we remove $e^{\sqrt{t}\Lambda}$ in the following Lemmas \ref{NSlem1} and \ref{NSlem2}, the results were essentially obtained in \cite{Wang04} (see also \cite{WaHuHaGu11}, Section 2.2.2). However, $e^{\sqrt{t}\Lambda}$ is of importance for the nonlinear estimates in the space $\widetilde{L}^{\gamma}(\mathbb{R}^{+}; e^{\sqrt{t}\Lambda} \dot{B}^{s}_{p,q})$, see Lemari\'{e}-Rieusset \cite{LeRi02} and \cite{BaBiTa12}.

\begin{lem}\label{NSlem1}
Let $1 < p < \infty$,  $1\leq q \leq \infty$,  $1\leq \gamma
\leq \infty$. Then
\begin{align} \label{2.1}
\|U_2(t) u_0\|_{\widetilde{L}^{\gamma}(\mathbb{R}^{+}; e^{\sqrt{t}\Lambda} \dot{B}^{s}_{p,
q})}\lesssim \|u_0 \|_{\dot{B}^{s-2/{\gamma}}_{p, q}}.
\end{align}
\end{lem}
{\it Proof.} In view of $e^{-t|\xi|^2/2 \ + \sqrt{t} |\xi|_1} \in M_p$ (see \cite{BaBiTa12}),  one easily sees that
\begin{align}
\|\Delta_j e^{\sqrt{t}\Lambda} U_2(t) u_0\|_{p} \le  \|\mathscr{F}^{-1}e^{-t|\xi|^2/2 \ + \sqrt{t} |\xi|_1} \|_1  \|\Delta_j U_2(t/2)
u_0\|_{p}  \lesssim   \|\Delta_j U_2(t/2)
u_0\|_{p}.   \label{NSfreloca}
\end{align}
We have known that (see \cite{Ch99} and \cite{WaHuHaGu11}, Section 2.2.2)
\begin{align}
\|\Delta_j U_2(t/2)f\|_{L^p_x} \lesssim e^{-ct2^{2j}}\|\Delta_j
f\|_{L^p_x},  \label{NSfreloc}
\end{align}
Hence, it follows from \eqref{NSfreloca} and \eqref{NSfreloc} that
\begin{align}
\|\Delta_j e^{\sqrt{t}\Lambda} U_2(t) u_0\|_{p}  \lesssim e^{-ct2^{2j}}  \|\Delta_j
u_0\|_{p}.   \label{NSfrelocb}
\end{align}
Taking $L_t^{\gamma}$ norm on inequality \eqref{NSfrelocb},
\begin{align}
\|\Delta_j e^{\sqrt{t}\Lambda} U_2(t) u_0\|_{L^\gamma_tL^p_x} \lesssim
2^{-2j/{\gamma}}\|\Delta_{j} u_0\|_{p}.\label{NSfreloc1}
\end{align}
Taking sequence $l^q$ norms in \eqref{NSfreloc1}, we have the result, as desired. $\hfill\Box$

\begin{lem}\label{NSlem2}
Let $1<p<\infty$, $1\leq   q\leq \infty$, $1\leq \gamma_1 \leq \gamma \le \infty$, $I= [0,T)$, $T\le \infty$. Then
\begin{align}
\|\mathscr{A}_2 f\|_{\widetilde{L}^{\gamma}(I;
e^{\sqrt{t}\Lambda} \dot{B}^{s}_{p,q})}\lesssim
\|f\|_{\widetilde{L}^{\gamma_1}(I;e^{\sqrt{t}\Lambda} \dot{B}^{s-2(1+1/\gamma -1/\gamma_1)}_{p,q})}.
\end{align}
In particular,
\begin{align}
\|\mathscr{A}_2 f\|_{\widetilde{L}^{\gamma}(I;
e^{\sqrt{t}\Lambda} \dot{B}^{s}_{p,q})}\lesssim
\|f\|_{\widetilde{L}^{1}(I;e^{\sqrt{t}\Lambda} \dot{B}^{s-2/\gamma}_{p,q})}. \label{2.7}
\end{align}
\end{lem}
{\it Proof.} In view of $e^{ (\sqrt{t}- \sqrt{\tau} - \sqrt{t-\tau}) |\xi|_1} \in M_p$ (see \cite{BaBiTa12}) and \eqref{NSfrelocb},
\begin{align} \label{2.8}
\|\Delta_j e^{\sqrt{t}\Lambda} \mathscr{A}_2 f\|_{p} & \lesssim \int_0^t
 \|\Delta_j e^{ (\sqrt{t}- \sqrt{\tau} )\Lambda + (t-\tau)\Delta} e^{  \sqrt{\tau} \Lambda} f(\tau)\|_{p}d\tau \nonumber\\
 & \lesssim \int_0^t
 \|\Delta_j e^{ \sqrt{t-\tau}\Lambda + (t-\tau)\Delta} e^{  \sqrt{\tau} \Lambda} f(\tau)\|_{p}d\tau \nonumber\\
  & \lesssim \int_0^t e^{-c(t-\tau)2^{2j}}
 \|\Delta_j  e^{\sqrt{\tau} \Lambda} f(\tau)\|_{p} d\tau.
\end{align}
Using Young's inequality, for $\gamma_1\le \gamma$,
\begin{align} \label{2.9}
\|\Delta_j e^{\sqrt{t}\Lambda} \mathscr{A}_2 f\|_{L^\gamma_{t\in I}L^p_x}
    \lesssim   2^{- 2j (1+1/\gamma -1/\gamma_1) }
 \|\Delta_j  e^{\sqrt{\tau} \Lambda} f(\tau)\|_{L^{\gamma_1}_{\tau \in I}L^p_x }.
\end{align}
Taking $l^q$ norms in both sides of \eqref{2.9}, we have the results, as desired. $\hfill\Box$

 For two
tempered distributions $f, g$, using the para-product decomposition \cite{CH, LeRi02}, one can decompose $fg$ as the
summation of the following two parts:
\begin{align*}
fg =\sum_{i,j}\Delta_i f \Delta_j g
   =\sum_{j}S_{j-1}f\Delta_j g+\sum_{i}S_ig \Delta_i f,
\end{align*}
By the support property in frequency spaces, one can easily check that
\begin{align*}
\Delta_i(S_{j}f\Delta_j g)=0, \,\,\textrm{for}\quad i>j+3.
\end{align*}
Therefore
\begin{align}
\Delta_i(fg)=\sum_{i\leq j+3} \Delta_i(S_{j-1}f\Delta_j g+S_j g
\Delta_j f).\label{para1}
\end{align}

\begin{lem}\label{NSlem3}
Let $1 < p   < \infty$,   $1\leq q, \gamma, \gamma_1, \gamma_2
\leq \infty$ with $1/\gamma=1/\gamma_1+1/\gamma_2$, $s, \varepsilon>0$. Then
\begin{align} \label{2.10}
\|fg\|_{\widetilde{L}^{\gamma}(I; e^{\sqrt{t}\Lambda} \dot{B}^{s}_{p,
q})} & \lesssim \|f\|_{\widetilde{L}^{\gamma_1}(I; e^{\sqrt{t}\Lambda} \dot{B}^{n/p-\varepsilon}_{p,
q})} \|g\|_{\widetilde{L}^{\gamma_2}(I; e^{\sqrt{t}\Lambda} \dot{B}^{s+\varepsilon}_{p,q})}  \nonumber\\
& \quad + \|g\|_{\widetilde{L}^{\gamma_1}(I; e^{\sqrt{t}\Lambda} \dot{B}^{n/p-\varepsilon}_{p,q})} \|f\|_{\widetilde{L}^{\gamma_2}(I; e^{\sqrt{t}\Lambda} \dot{B}^{s+\varepsilon}_{p,q})}.
\end{align}
Moreover, if $\varepsilon=0$, then \eqref{2.10} also holds for $q=1.$
\end{lem}
{\it Proof.} Denote $F= e^{\sqrt{t}\Lambda} f, \ G= e^{\sqrt{t}\Lambda} g$.  In view of \eqref{para1},
\begin{align}
e^{\sqrt{t}\Lambda} \Delta_i(fg) & =\sum_{i-3\leq j } e^{\sqrt{t}\Lambda} \Delta_i(S_{j-1}f\Delta_j g+S_j g
\Delta_j f) \nonumber\\
& =\sum_{i-3\leq j } \ \sum_{k\le j-1} e^{\sqrt{t}\Lambda} \Delta_i(\Delta_{k}e^{-\sqrt{t}\Lambda} F  \ \Delta_j e^{-\sqrt{t}\Lambda} G) \nonumber\\ & \quad +  \sum_{i-3\leq j } \ \sum_{k\le j} e^{\sqrt{t}\Lambda} \Delta_i(\Delta_{k}e^{-\sqrt{t}\Lambda} G  \ \Delta_j e^{-\sqrt{t}\Lambda} F)
\label{2.11}
\end{align}
Now we use an idea as in \cite{LeRi02}, Page 253 and  \cite{BaBiTa12}, to consider the bilinear form
\begin{align}
 B_t(u, v) & = e^{\sqrt{t}\Lambda}(e^{-\sqrt{t}\Lambda}u \
e^{-\sqrt{t}\Lambda}v)  \nonumber \\
 & =
\int_{\mathbb{R}^n}\int_{\mathbb{R}^n}e^{ix
\xi}e^{\sqrt{t}(|\xi|_1-|\xi-\eta|_1-|\eta|_1)}\hat{u}(\xi-\eta)\hat{v}(\eta)d{\eta}d{\xi}. \label{2.12}
\end{align}
Denote for $\lambda=(\lambda_1,...,\lambda_n), \mu=(\mu_1,...,\mu_n), \nu=(\nu_1,...,\nu_n)$, $ \lambda_i,  \mu_i , \nu_i \in \{1, -1\} $,
\begin{align}
 D_\lambda & = \{\eta: \ \lambda_i \eta_i \ge 0, \ i=1,...,n\}, \nonumber\\
    D_\mu & = \{\xi-\eta: \ \mu_i (\xi_i-\eta_i) \ge 0, \ i=1,...,n\},\nonumber\\
       \ D_\nu & = \{\xi: \ \nu_i  \xi_i  \ge 0, \ i=1,...,n\}. \nonumber
\end{align}
We denote by $\chi_D$ the characteristic function on $D$. Then we can rewrite $B_t(u,v)$ as
\begin{align}
  \sum_{\lambda_i, \mu_j, \nu_k \in \{1,-1\}}
\int_{\mathbb{R}^n}\int_{\mathbb{R}^n} e^{ix
\xi}\chi_{D_\nu}(\xi)e^{\sqrt{t}(|\xi|_1-|\xi-\eta|_1-|\eta|_1)}\chi_{D_\mu}(\xi-\eta)\hat{u}(\xi-\eta)\chi_{D_\lambda}(\eta)\hat{v}(\eta)d{\eta}d{\xi}.
\end{align}
Noticing that for $\eta\in D_\lambda$, $\xi-\eta\in D_\mu$ and $\xi \in D_\mu$,  $e^{\sqrt{t}(|\xi_i|-|\xi_i-\eta_i|-|\eta_i|)}$ must belong to the following set
\begin{align}
\mathfrak{M} =\{ \ 1, \ e^{-2 \sqrt{t}|\xi_i|}, \  e^{- 2 \sqrt{t} |\xi_i-\eta_i|} ,  e^{- 2 \sqrt{t} |\eta_i|} \}, \quad i=1,...,n.
\end{align}
For $1<p<\infty$,  $\chi_{D_\lambda} \in M_p$, $m\in M_p$ for any $m\in \mathfrak{M}$, it follows from the algebra property of $M_p$ that
\begin{align}
\| B_t(u, v)\|_p \lesssim \|u\|_{p_1} \|v\|_{p_2}, \ \ 1/p_1 +1/p_2=1. \label{2.15}
\end{align}
By \eqref{2.11}, \eqref{2.12} and \eqref{2.15}, then using Bernstein's inequality, we have
\begin{align}
2^{si}\|e^{\sqrt{t}\Lambda} \Delta_i(fg)\|_{L^\gamma_{t\in I} L^p_x }
&  \lesssim \sum_{i-3\leq j } \ \sum_{k\le j-1} 2^{s(i- j)} \| \Delta_{k}  F \|_{L^{\gamma_1}_{t\in I} L^\infty_x } 2^{sj}  \| \Delta_j   G\|_{L^{\gamma_2}_{t\in I} L^{p}_x }  \nonumber\\
 & \quad +  \sum_{i-3\leq j } \ \sum_{k\le j} 2^{s(i-j)} \|\Delta_{k} G\|_{L^{\gamma_1}_{t\in I} L^\infty_x }  2^{sj} \| \Delta_j   F\|_{L^{\gamma_2}_{t\in I} L^{p}_x } \nonumber\\
& \lesssim \| F \|_{\widetilde{L}^{\gamma_1} (I; \dot{B}^{n/p-\varepsilon}_{p,q}) } \sum_{i-3\leq j } \ 2^{s(i- j)}   2^{(s+\varepsilon)j}  \| \Delta_j   G\|_{L^{\gamma_2}_{t\in I} L^{p}_x }  \nonumber\\
 & \quad +  \| G \|_{\widetilde{L}^{\gamma_1} (I; \dot{B}^{n/p-\varepsilon}_{p,q}) }\sum_{i-3\leq j } \ 2^{s(i-j)}   2^{(s+\varepsilon)j} \| \Delta_j   F\|_{L^{\gamma_2}_{t\in I} L^{p}_x }.
\label{2.16}
\end{align}
Taking sequence $l^q$ norm over all $i\in \mathbb{Z}$ in \eqref{2.16} and using Young's inequality, we immediately have the result, as desired. $\hfill\Box$

\begin{cor}\label{NSlem4}
Let $1< p   < \infty$, $1\le q\le \infty$. Then
\begin{align} \label{2.17}
\|fg\|_{\widetilde{L}^{1}(I; e^{\sqrt{t}\Lambda} \dot{B}^{n/p}_{p,
q})} & \lesssim \|f\|_{\widetilde{L}^{3}(I; e^{\sqrt{t}\Lambda} \dot{B}^{n/p-1/3}_{p,
q})} \|g\|_{\widetilde{L}^{3/2}(I; e^{\sqrt{t}\Lambda} \dot{B}^{n/p+1/3}_{p,q})} \nonumber\\
& \quad \ + \|g\|_{\widetilde{L}^{3}(I; e^{\sqrt{t}\Lambda} \dot{B}^{n/p-1/3}_{p,
q})} \|f\|_{\widetilde{L}^{3/2}(I; e^{\sqrt{t}\Lambda} \dot{B}^{n/p+1/3}_{p,q})}
\end{align}
\end{cor}

{\it Proof of Theorem \ref{NSTh1a}.} Now let us consider the map
\begin{align} \label{2.18}
\mathscr{T}: u(t) \to U_2(t) u_0 + \mathscr{A}_2 \mathbb{P}{\rm div} (u\otimes u)
\end{align}
in the metric space ($I=[0,T]$)
\begin{align*}
& \mathscr{D} = \{u: \|u\|_{\widetilde{L}^{3}(I; e^{\sqrt{t}\Lambda} \dot{B}^{n/p-1/3}_{p,q}) \cap \widetilde{L}^{3/2}(I; e^{\sqrt{t}\Lambda} \dot{B}^{n/p+1/3}_{p,q})} \le \delta\},  \\
& \quad d(u,v) =\|u-v\|_{\widetilde{L}^{3}(I; e^{\sqrt{t}\Lambda} \dot{B}^{n/p-1/3}_{p,q}) \cap \widetilde{L}^{3/2}(I; e^{\sqrt{t}\Lambda} \dot{B}^{n/p+1/3}_{p,q})}.
\end{align*}
By   \eqref{2.7} and Corollary \ref{NSlem4}, for any $u, v\in \mathscr{D}$,
\begin{align} \label{2.20}
\|\mathscr{T} & u\|_{\widetilde{L}^{3}(I; e^{\sqrt{t}\Lambda} \dot{B}^{n/p-1/3}_{p,q}) \cap \widetilde{L}^{3/2}(I; e^{\sqrt{t}\Lambda} \dot{B}^{n/p+1/3}_{p,q})} \nonumber\\
  & \lesssim  \|U_2(t) u_0\|_{\widetilde{L}^{3}(I; e^{\sqrt{t}\Lambda} \dot{B}^{n/p-1/3}_{p,q}) \cap \widetilde{L}^{3/2}(I; e^{\sqrt{t}\Lambda} \dot{B}^{n/p+1/3}_{p,q}) } \nonumber\\
  & \quad +  \|u\|^2_{\widetilde{L}^{3}(I; e^{\sqrt{t}\Lambda} \dot{B}^{n/p-1/3}_{p,q}) \cap \widetilde{L}^{3/2}(I; e^{\sqrt{t}\Lambda} \dot{B}^{n/p+1/3}_{p,q})} \nonumber\\
   &  \lesssim   \|U_2(t) u_0\|_{\widetilde{L}^{3}(I; e^{\sqrt{t}\Lambda} \dot{B}^{n/p-1/3}_{p,q}) \cap \widetilde{L}^{3/2}(I; e^{\sqrt{t}\Lambda} \dot{B}^{n/p+1/3}_{p,q}) }  +  \delta^2, \\
& d(\mathscr{T} u, \  \mathscr{T}v ) \lesssim  \delta d( u, v) . \label{2.21}
\end{align}
There exists a $T>0$ such that
\begin{align} \label{2.22}
      \|U_2(t) u_0\|_{\widetilde{L}^{3}(I; e^{\sqrt{t}\Lambda} \dot{B}^{n/p-1/3}_{p,q}) \cap \widetilde{L}^{3/2}(I; e^{\sqrt{t}\Lambda} \dot{B}^{n/p+1/3}_{p,q})} \le   \delta /2.
\end{align}
Indeed, we can choose a sufficiently large $J\in \mathbb{N}$ such that
$$\left(\sum_{|j|\ge J} 2^{q(n/p-1)j}\|\triangle_j u_0\|^q_{p}\right)^{1/q} \le \delta/4 .$$
It follows from \eqref{NSfreloc1} that
\begin{align*}
   \left(\sum_{|j|\ge J}  2^{q(n/p-1/3)j} \|\triangle_j e^{\sqrt{t}\Lambda} U_2(t) u_0\|^q_{ L^{3}_{t\in I}L^p_x} \right)^{1/q} \le  \left(\sum_{|j|\ge J} 2^{q(n/p-1)j}\|\triangle_j u_0\|^q_{p}\right)^{1/q} \le \delta/4 .
\end{align*}
On the other hand, one can choose $|I|$ small enough such that
\begin{align*}
   \left(\sum_{|j|\le J}  2^{q(n/p-1/3)j} \|\triangle_j e^{\sqrt{t}\Lambda} U_2(t) u_0\|^q_{ L^{3}_{t\in I}L^p_x} \right)^{1/q}    \le \delta/4.
\end{align*}
So, we have \eqref{2.22}.
Applying the standard  argument, we can show that $\mathscr{T}$ is a contraction mapping from $(\mathscr{D},d) $ into itself.  So, there is a $u\in \mathscr{D}$ satisfying $\mathscr{T} u=u$.  Moreover, by Lemmas \ref{NSlem1} and \ref{NSlem2},
\begin{align} \label{2.23}
   \|  u\|_{\widetilde{L}^{\infty}(I; e^{\sqrt{t}\Lambda} \dot{B}^{n/p-1}_{p,q})} & \lesssim  \|u_0\|_{\dot{B}^{n/p-1}_{p,q}}  +\delta^2.
\end{align}
The solution can be extended step by step and finally we have a maximal time $T_{\max}$ verifying
$$
u \in  \widetilde{L}^{3}_{\rm loc}(0,T_{\max}; e^{\sqrt{t}\Lambda} \dot{B}^{n/p-1/3}_{p,q}) \cap \widetilde{L}^{3/2}_{\rm loc}(0,T_{\max}; e^{\sqrt{t}\Lambda} \dot{B}^{n/p+1/3}_{p,q}) \cap  \widetilde{L}^{\infty}([0,T_{\max}); e^{\sqrt{t}\Lambda} \dot{B}^{n/p-1}_{p, q} ).
$$
If $T_{\max}<\infty$ and $\|u \|_{\widetilde{L}^{3}(0,T_{\max}; e^{\sqrt{t}\Lambda} \dot{B}^{n/p-1/3}_{p,q}) \cap \widetilde{L}^{3/2}(0,T_{\max}; e^{\sqrt{t}\Lambda} \dot{B}^{n/p+1/3}_{p,q})}<\infty$, we  claim that the solution can be extended beyond $T_{\max}$. Indeed, let us consider the integral equation
\begin{align} \label{2.24}
u(t) =   U_2(t-T) u (T) + \int^t_T U_2(t-\tau) \mathbb{P}{\rm div} (u\otimes u)(\tau) d\tau,
\end{align}
we see that
\begin{align} \label{2.25}
 \|U_2&(t-T) u (T)\|_{\widetilde{L}^{3}(T,T_{\max};  e^{\sqrt{t}\Lambda}\dot{B}^{n/p-1/3}_{p, q}) \cap \widetilde{L}^{3/2}(T,T_{\max};  e^{\sqrt{t}\Lambda}\dot{B}^{n/p+1/3}_{p, q})} \nonumber\\
  & \le  \|u\|_{\widetilde{L}^{3}(T,T_{\max};  e^{\sqrt{t}\Lambda}\dot{B}^{n/p-1/3}_{p, q}) \cap \widetilde{L}^{3/2}(T,T_{\max};  e^{\sqrt{t}\Lambda}\dot{B}^{n/p+1/3}_{p, q})} \nonumber\\
   & \quad \  + C \|u\|^2_{\widetilde{L}^{3}(T,T_{\max};  e^{\sqrt{t}\Lambda}\dot{B}^{n/p-1/3}_{p, q}) \cap \widetilde{L}^{3/2}(T,T_{\max};  e^{\sqrt{t}\Lambda}\dot{B}^{n/p+1/3}_{p, q})}
   \le \delta/2
\end{align}
if $T$ is sufficiently close to $T_{\max}$. \eqref{2.25} is analogous to \eqref{2.22}, which implies that the solution exists on $[T,T_{\max}]$. A contradiction with $T_{\max}$ is maximal. Moreover, if $\|u_0\|_{\dot{B}^{n/p-1}_{p,q}} \ll 1$, we can directly choose $T=\infty$ in \eqref{2.20} and \eqref{2.21}. $\hfill\Box$

\section{Gevrey Class and $E^s_{p, q}$} \label{sectGevExp}

In \cite{Wa06}, it was shown that $G_{1, 2} = \bigcup_{s>0} {E}^{s}_{2, q}.$ In this paper we generalize this equality to all  $p>0$.

\begin{prop}\label{prop4.1}
{\rm ({ \cite{Tr},  Nikol'skij's inequality})} \it Let $\Omega\subset
\mathbb{R}^n$ be a compact set,  $0<r\le \infty.$ Let us denote
$\sigma_r=n( 1/(r\wedge 1)-1/2)$ and assume that $s>\sigma_r$. Then
there exists a constant $C>0$ such that
\begin{align*}
\|\mathscr{F}^{-1}\varphi \mathscr{F} f\|_r\le
C\|\varphi\|_{H^s}\|f\|_r
\end{align*}
holds for all $f\in L^r_\Omega:=\{f\in L^p: \; {\rm supp}
\widehat{f} \subset \Omega\}$ and $\varphi\in H^s$. Moreover,  if
$r\ge 1$,  then the above inequality holds for all $f\in L^r$.
\end{prop}

\begin{prop}\label{prop4.2}
Let $ 0<  p ,  q \leq \infty$. Then
\begin{align}
 G_{1, p} = \bigcup_{s>0} {E}^{s}_{p, q}.
\end{align}
\end{prop}

{\it Proof.} We have
\begin{align}
 \|\partial^\alpha f \|_p \le \sum_{k\in \mathbb{Z}^n} \|\Box_k \partial^\alpha f \|_p.
\end{align}
We easily see that
\begin{align}
  \|\Box_k \partial^\alpha f \|_p \le \sum_{|l|_\infty\le 1} \|\sigma_{k+l} \xi^\alpha \|_{M_p} \|\Box_k f\|_p.
\end{align}
Since $M_p$ is translation invariant,  in view of Nikol'skij's inequality we have
\begin{align}
 \|\sigma_{k} \xi^\alpha \|_{M_p} = \|\sigma \cdot (\xi+k)^\alpha \|_{M_p} \lesssim  \|\sigma \cdot (\xi+k)^\alpha \|_{H^L},  \ \ L>(n/(1\wedge p) -1/2).
\end{align}
It is easy to calculate that
\begin{align}
  \|(\xi+k)^\alpha \sigma   \|_{H^L} & \lesssim \|(\xi+k)^\alpha \sigma\|_2 +  \sum^n_{i=1} \|\partial^L_{x_i}((\xi +k)^\alpha \sigma)\|_2 \nonumber\\
& \lesssim \prod^n_{i=1} \langle k_i\rangle^{\alpha_i} + \sum^n_{i=1} \sum^L_{l_i=1}  \prod^{l_i-1}_{m=0} (\alpha_i -m)   \prod^n_{j=1,  j\neq i} \langle k_j\rangle^{ \alpha_j} \langle k_i\rangle^{ \alpha_i-l_i}.
\end{align}
It follows that
\begin{align}
 \|\partial^\alpha f \|_p
 & \lesssim     \sum_{k\in \mathbb{Z}^n}   \prod^n_{i=1} \langle k_i\rangle^{\alpha_i}   \|\Box_k   f \|_p \nonumber\\
  & + \sum^n_{i=1} \sum^L_{l_i=1} \ \sum_{k\in \mathbb{Z}^n} \prod^{l_i-1}_{m=0} (\alpha_i -m)   \prod^n_{j=1,  j\neq i} \langle k_j\rangle^{ \alpha_j} \langle k_i\rangle^{ \alpha_i-l_i}  \|\Box_k   f \|_p. \nonumber
\end{align}
If we can show that
\begin{align} \label{3.8}
  \frac{s^{|\alpha|}} {\alpha !} \frac{\prod^{l_i-1}_{m=0} (\alpha_i -m)   \prod^n_{j=1,  j\neq i} \langle k_j\rangle^{ \alpha_j} \langle k_i\rangle^{ \alpha_i-l_i} }{ e^{2s(|k_1|+...+|k_n|)}} \lesssim 1,
\end{align}
then together with
 $ E^{s+ \varepsilon}_{p, q} \subset E^{s}_{p, 1}$ for all $ s,  \varepsilon>0$ (cf. \cite{Wa06}),   we immediately have
 \begin{align} \label{3.9}
 \|\partial^\alpha f \|_p
 \lesssim   \frac{\alpha !} {s^{|\alpha|}}  \sum_{k\in \mathbb{Z}^n}  e^{2s(|k_1|+...+|k_n|)}     \|\Box_k   f(x)\|_p \lesssim   \frac{\alpha !} {s^{|\alpha|}}       \|f\|_{E^{cs}_{p, 1}} \lesssim   \frac{\alpha !} {s^{|\alpha|}}       \|f\|_{E^{\tilde{c}s}_{p, q}}
\end{align}
for any $ \tilde{c}>c= 2n \mathrm{log}_2 e$.   Noticing that $s>0$ is arbitrary,   it follows from \eqref{3.9} that
\begin{align}
\bigcup_{s>0} {E}^{s}_{p, q} \subset  G_{1, p} .
\end{align}

Now we show \eqref{3.8}.  Using Taylor's expansion,  we see that
\begin{align}
 e^{2s |k_i| } \ge \frac{(2s |k_i|)^{ \alpha_i -l_i  }}{(\alpha_i -l_i) !},  \quad e^{2s |k_j| } \ge \frac{(2s |k_j|)^{ \alpha_j }}{ \alpha_j !},
 \end{align}
 from which we see that \eqref{3.8} holds.

Next,  we show that   $ G_{1, p}  \subset \bigcup_{s>0} E^s_{p,  \infty}$.  We have
\begin{align} \label{3.12}
 \|f\|_{E^{s\mathrm{log}_2 e}_{p,  \infty}} & = \sup_k e^{s|k|} \|\Box_k f\|_p \nonumber\\
 & \lesssim   \sup_k  \sum^\infty _{m=0}  \frac{  (ns)^{ m} }{ m !} (|k_1 |^m+...+ |k_n|^m) \|\Box_k f\|_p
\end{align}
If $|k_i| \le 10$,  we see that
\begin{align} \label{3.13}
    \sum^\infty _{m=0}  \frac{ (ns)^{ m} }{ m !}   |k_i|^m  \|\Box_k f\|_p  \lesssim   \sum^\infty _{m=0}  \frac{  (10ns)^{ m} }{ m !}      \| f\|_p \le M e^{10ns}.
\end{align}
 If $|k_i| > 10$,  in view of Nikol'skij's inequality,  for some $L>(n/(1\wedge p) -1/2)$,
\begin{align} \label{3.14}
    |k_i|^m  \|\Box_k f\|_p & \lesssim  |k_i|^m \sum_{|l|_\infty\le 1} \|\sigma_{k+l} |\xi_i|^{-m} \|_{M_p} \|\Box_k \partial^{m}_{x_i} f\|_p \nonumber\\
      & \lesssim   2^m (1 + m \langle k_i\rangle^{-1}+...+  m(m+1)...(m+L-1) \langle k_i\rangle^{-L})  \|  \partial^{m}_{x_i} f\|_p \nonumber\\
       & \lesssim   2^m  (1+ (m+L)^{L})  \|  \partial^{m}_{x_i} f\|_p .
\end{align}
By \eqref{3.12} and \eqref{3.14},  if $f\in G_{1, p}$ and $|k_i| > 10$,   then for any $s<\rho /2n$,
\begin{align} \label{3.15}
    \sum^\infty _{m=0}  \frac{  (ns)^{ m} }{ m !}  |k_i |^m  \|\Box_k f\|_p  \lesssim \sum^\infty _{m=0} \frac{(2ns)^{m}}{\rho^m } (1+(m+L)^{L})    \lesssim 1.
\end{align}
By \eqref{3.13} and \eqref{3.15},  we have $f\in E^s_{p, \infty}$. It follows that $ G_{1, p}  \subset \bigcup_{s>0} E^s_{p,  \infty}$.  Noticing that $E^s_{p,  \infty} \subset E^{s/2}_{p,  q}$,  we immediately have $ G_{1, p}  \subset \bigcup_{s>0} E^s_{p,  q}$. $\hfill\Box$

\section{Analyticity of NS: Proof of Theorem \ref{NSTh1b}} \label{sectNS}

In Section \ref{sectNSa}, the nonlinear mapping estimate is invalid for the case $p=\infty$,  we cannot use the same way as in Section \ref{sectNSa} to obtain the analyticity of the solutions if $p=\infty$. Our idea is to use the frequency-uniform decomposition techniques to
handle the case $p=\infty$, where the corresponding space for the initial data should be modulation spaces. The local well posedness of NS in $M^{-1}_{\infty,1}$ was established in \cite{Iw10}.

Let $1\le p \le \infty$.  We show that there exists $c>0$ (say $0<c \le 2^{-10}$) such that
\begin{align} \label{4.1}
\|\Box_k    U_2(t) f\|_{p}  \lesssim  e^{-2ct|k|^2} \|\Box_k f \|_{p}
\end{align}
holds for all $f\in L^p$ and $k\in \mathbb{Z}^n$. Let $\tilde{\sigma}:\mathbb{R}^n \to [0,1]$ be a smooth cut-off function satisfying $\tilde{\sigma} (\xi)=1$ for $|\xi|_\infty \le 3/4$ and $\tilde{\sigma} (\xi)=0$ for $|\xi|_\infty >7/8$. Let us observe  that
\begin{align}
\|\Box_k    U_2(t) f\|_{p} & =    \|\mathscr{F}^{-1} \tilde{\sigma} (\xi-k) e^{- t|\xi|^2} \  \widehat{\Box_k f}\|_{p} \nonumber\\
& \lesssim    \|  \tilde{\sigma} (\xi-k) e^{- t|\xi|^2}\|_{M_p} \|\Box_k f \|_{p}.
\end{align}
In view of Nikol'skij's inequality,
\begin{align}
  \| \tilde{\sigma} (\xi-k) e^{- t|\xi|^2}\|_{M_p}   \lesssim \| \tilde{\sigma}  e^{- t|\xi+k|^2}\|_{H^L},  \quad L>n(1/(1\wedge p)-1/2).
\end{align}
Since for $|k| \ge 1$,
\begin{align}
 & \partial^L_{\xi_i} (\tilde{\sigma}  e^{P(\xi)})  =    \sum_{\beta+\gamma =L} C_{L\beta\gamma} \partial^\beta_{\xi_i} \tilde{\sigma} \ \partial^\gamma_{\xi_i}  e^{P(\xi)}, \label{4.4}\\
 &   |\partial^\gamma_{\xi_i}  e^{P(\xi)}| \lesssim  e^{P(\xi)} \sum_{\substack{\alpha_1+...+\alpha_q=\gamma \\ 1\le q \le \gamma}} |\partial^{\alpha_1}_{\xi_i} P(\xi) ... \partial^{\alpha_q}_{\xi_i} P(\xi) |, \label{4.5}
\end{align}
we immediately have for $|k|  \ge 1$,  $| \partial^\gamma_{\xi_i}   e^{- t|\xi+k|^2} |  \lesssim  e^{- t|\xi+k|^2/2}$, so,
\begin{align} \label{4.6}
  \| \tilde{\sigma} (\xi-k) e^{- t|\xi|^2}\|_{M_p}   \lesssim  e^{- t|k|^2/ 2^{7} }.
\end{align}
We easily see that  \eqref{4.6} also holds for $k=0$. Hence,   we have shown \eqref{4.1}.

Now let $0<t_0 \le \infty$.  We consider the estimate of $  \mathscr{A}_2  f  $. By \eqref{4.1},
we have
\begin{align} \label{4.8}
 \|\Box_k  \mathscr{A}_2  f\|_{L^p_{x}} \lesssim \int^t_0   2^{-2c(t-\tau)|k|^2} \|\Box_k f (\tau)\|_{L^p_{x}} d\tau.
\end{align}
It follows from \eqref{4.8} that
\begin{align} \label{4.9}
 2^{ct|k|} \|\Box_k  \mathscr{A}_2  f\|_{L^p_{x}} \lesssim \int^t_0  2^{- c(t-\tau)|k|^2} 2^{ c\tau |k|}\|\Box_k f (\tau)\|_{L^p_{x}} d\tau.
\end{align}
By Young's inequality,  from \eqref{4.9} we have
\begin{align} \label{4.10}
 \| 2^{ct|k|} \Box_k  \mathscr{A}_2  f\|_{L^{2}_{ t\in [0, t_0]}L^p_{x}} \lesssim  \langle k\rangle ^{- 1} \|2^{ c\tau |k|}\Box_k f (\tau)\|_{L^{1}_{ t\in [0, t_0]}L^p_{x}}.
\end{align}
Applying $\|\nabla \Box_k     f\|_{L^p_{x}}  \lesssim \langle k\rangle \|\Box_k f\|_{L^p_{x}}$,  we have
\begin{align} \label{4.11}
 \| 2^{ct|k|} \nabla \Box_k  \mathscr{A}_2 f\|_{L^{2}_{t\in [0, t_0]}L^p_{x}} \lesssim    \|2^{ ct |k|}\Box_k f\|_{L^{1}_{t\in [0, t_0]}L^p_{x}}.
\end{align}
For convenience,  we denote for any $\mathbb{A} \subset \mathbb{Z}^n$,
\begin{align} \label{4.12}
\|f\|_{\widetilde{L}^{{{\tilde{q}}}}(I, \mathbb{A};E^{s}_{p, q})} = \left(\sum_{k\in  \mathbb{A}} \|2^{s|k|}\Box_k f\|^q_{L^{{\tilde{q}}}_{t\in I}L^p_{x}}\right)^{1/q}.
\end{align}
So,  taking the sequence $l^q$ norm in both sides of \eqref{4.11},  we have

\begin{prop} \label{prop5.2}
Let $1\le p\le \infty$. There exists a constant $c>0$ ($0<c \le 2^{-10}$) such that for $0< t_0 < \infty$ and $I=[0, t_0]$,
\begin{align} \label{4.12}
\| U_2(t) u_0\|_{\widetilde{L}^{2}(I, \mathbb{A};E^{ct}_{p, 1}) } & \lesssim  \sum_{k\in \mathbb{A}} \langle k\rangle^{-1} \|\Box_k u_0 \|_{p}, \quad \mathbb{A} \subset \mathbb{Z}^n \setminus\{0\},\\
\|\Box_0    U_2(t) u_0\|_{L^{2}_{t\in I} L^p_x }  & \lesssim   |I|^{1/2} \|\Box_k u_0 \|_{p},  \label{4.13}\\
\|\nabla   \mathscr{A}_2 f\|_{\widetilde{L}^{2}(I, \mathbb{A};E^{ct}_{p, 1})}  & \lesssim  \|f\|_{\widetilde{L}^{1}(I, \mathbb{A};E^{ct}_{p,1})}, \quad \mathbb{A} \subset \mathbb{Z}^n \setminus\{0\},\\
\|\Box_0   \mathscr{A}_2 f\|_{ L^{2}_{t\in I} L^p_x}  & \lesssim |I|^{1/2}  \|\Box_0 f\|_{ L^{1}_{t\in I} L^p_x },
\end{align}
  holds for all $u_0\in M^{-1}_{p,1}$ and  $f\in  \widetilde{L}^{ 1}(I;E^{ct}_{p, 1})$.
  \end{prop}

Now we consider a nonlinear mapping estimate in  $\widetilde{L}^{{{\tilde{q}}}}(I;E^{s(t)}_{p, q})$.
\begin{prop} \label{prop5.3}
Let $1\le p,  p_1,  p_2,  q,  q_1,  q_2,  \tilde{q},  \tilde{q}_1,  \tilde{q}_2 \le \infty$.  $1/p=1/p_1+1/p_2,  \ 1/\tilde{q}=1/ \tilde{q}_1 + 1/\tilde{q}_2$ and $1/q=1/q_1+1/q_2-1$, $s(\cdot):\mathbb{R}_+ \to \mathbb{R}_+$ and $I\subset \mathbb{R}_+$ with $\sup_{t\in I} s(t)<\infty$. Then we have
\begin{align} \label{4.14}
\|fg\|_{\widetilde{L}^{{{\tilde{q}}}}(I;E^{s(t)}_{p, q})}   \lesssim \sup_{t\in I} 2^{4s(t)} \|f\|_{\widetilde{L}^{{{\tilde{q}_1}}}(I;E^{s(t)}_{p_1, q_1})} \|g\|_{\widetilde{L}^{{{\tilde{q}_2}}}(I;E^{s(t)}_{p_2, q_2})},
\end{align}
 \end{prop}

{ \it Proof. } By definition,  we have
\begin{align} \label{4.15}
\|fg\|_{\widetilde{L}^{{{\tilde{q}}}}(I;E^{s(t)}_{p, q})} = \left(\sum_{k\in\mathbb{Z}^n} \|2^{s(t) |k|}\Box_k (fg)\|^q_{L^{{\tilde{q}}}_{t\in I}L^p_{x}}\right)^{1/q}.
\end{align}
Using the fact
\begin{align} \label{4.16}
 \Box_k (fg)= \sum_{i,  j\in {\Bbb Z}^n} \Box_k (\Box_i f \;
\Box_j g),  \ \  \Box_k (\Box_i f \;
\Box_j g) =0 \ \ \mathrm{if} \ \ |k-i-j| >4,
\end{align}
by H\"older's inequality, we have
\begin{align} \label{4.17}
  \|2^{s(t)|k|}\Box_k (fg)\|_{L^{{\tilde{q}}}_{t\in I}L^p_{x}} & \le \sum_{i,  j\in {\Bbb Z}^n}   \|2^{s(t)|k|}\Box_k (\Box_i f \;
\Box_j g) \|_{L^{{\tilde{q}}}_{t\in I}L^p_{x}} \chi_{|k-i-j|\le 4} \nonumber\\
& \lesssim  \sum_{|l|_\infty \le 4 }  \sum_{i \in {\Bbb Z}^n}  \|2^{s(t) (|i|+|k-i|)} \Box_i f \Box_{k-i-l} g
  \|_{L^{{\tilde{q} }}_{t\in I}L^{p }_{x}}   \nonumber\\
& \lesssim  \sum_{|l|_\infty \le 4 }  \sum_{i \in {\Bbb Z}^n} \| 2^{s(t)|i|}\Box_i f
  \|_{L^{{\tilde{q}_1}}_{t\in I}L^{p_1}_{x}}\|2^{s(t)|k-i|}\Box_{k-i-l} g
  \|_{L^{{\tilde{q}_2}}_{t\in I}L^{p_2}_{x}} \nonumber\\
& \lesssim \sup_{t\in I} 2^{4s(t)} \sum_{|l|_\infty \le 4 }  \sum_{i \in {\Bbb Z}^n} \| 2^{s(t)|i|}\Box_i f
  \|_{L^{{\tilde{q}_1}}_{t\in I}L^{p_1}_{x}} \nonumber\\
  & \quad\quad \times \|2^{s(t)|k-i-l|}\Box_{k-i-l} g
  \|_{L^{{\tilde{q}_2}}_{t\in I}L^{p_2}_{x}}.
\end{align}
Combining \eqref{4.15} and \eqref{4.17} and applying Young's inequality,  we immediately have the result,  as desired. $\hfill\Box$

Now let us consider the map
\begin{align} \label{4.18}
\mathscr{T}: u(t) \to U_2(t) u_0 + \mathscr{A}_2 \mathbb{P}{\rm div} (u\otimes u).
\end{align}
 Let $c=2^{-10}$, $I=[0,t_0]$. For any $\delta>0$, one  can choose $J:=J(u_0)>0$ satisfying  $\sum_{|k|\ge J} \|\Box_k u_0\|_{p} \le \delta/4C$. By Proposition \ref{prop5.2},
\begin{align} \label{4.19}
\|U_2(t) u_0\|_{\widetilde{L}^{2}(I, \{k: |k|\ge J\};E^{ct}_{p, 1})}   \le C \sum_{|k|\ge J} \|\Box_k u_0\|_{p} \le \delta/4.
\end{align}
On the other hand, one can choose $I$ satisfying $ |I|^{1/2} \langle J\rangle \le \delta/4C \|u_0\|_{M^{-1}_{p,1}}$, so,
\begin{align} \label{4.19a}
\|U_2(t) u_0\|_{\widetilde{L}^{2}(I, \{k: |k|< J\};E^{ct}_{p, 1})} &  \le C \sum_{|k|< J}  \left\|\| 2^{-ct|k|^2} \Box_k u_0\|_{p}
\right\|_{L^2_{t\in I}} \nonumber\\
&  \le C |I|^{1/2} \langle J\rangle \sum_{|k|< J} \langle k\rangle^{-1}  \|\Box_k u_0\|_{p} \le \delta/4.
\end{align}
Hence, in view of \eqref{4.19} and \eqref{4.19a},
\begin{align} \label{4.19aa}
\|U_2(t) u_0\|_{\widetilde{L}^{2}(I;E^{ct}_{p, 1})}   \le   \delta/2.
\end{align}
Denote $\widetilde{\Box}_0 = \mathscr{F}^{-1}\tilde{\sigma}\mathscr{F}$.  Noticing that $\widetilde{\Box}_0\mathbb{P}{\rm div}: L^p\to L^p$ for all $p\in [1,\infty]$,  we have from Propositions \ref{prop5.2} and \ref{prop5.3} that
\begin{align} \label{4.20}
\|\mathscr{A}_2 \mathbb{P}{\rm div} (u\otimes u)\|_{\widetilde{L}^{2}(I;E^{ct}_{p,1})}  &  \lesssim (1+|I|^{1/2}) \|u\otimes u\|_{\widetilde{L}^{1} (I;E^{ct}_{p,1})} \nonumber\\
 &  \lesssim  (1+|I|^{1/2})2^{4ct_0} \|u\|_{\widetilde{L}^{2} (I;E^{ct}_{p,1})}\|u\|_{\widetilde{L}^{2} (I;E^{ct}_{\infty,1})} \nonumber\\
 &  \lesssim (1+|I|^{1/2}) 2^{4ct_0} \|u\|^2_{\widetilde{L}^{2} (I;E^{ct}_{p,1})}.
\end{align}
We can assume that $|I|\le 1$.  Hence, collecting \eqref{4.19aa} and \eqref{4.20}, we have
\begin{align} \label{4.21}
\|\mathscr{T} u\|_{\widetilde{L}^{2}(I;E^{ct}_{p,1})}  &  \le \delta/2  + C \|u\|^2_{\widetilde{L}^{2} (I;E^{ct}_{p,1})}.
\end{align}
Now we can fix $\delta$ verifying $C\delta \le 1/4.$  Put
\begin{align} \label{4.21a}
 \mathscr{D} =\{u\in  \widetilde{L}^{2}(I;E^{ct}_{p,1})  : \  \|u\|_{\widetilde{L}^{2} (I;E^{ct}_{p,1})} \le \delta\}
\end{align}
We have $ \mathscr{T} u \in  \mathscr{D}$ if  $ u \in  \mathscr{D}$ and
\begin{align} \label{4.23}
  \|\mathscr{T} u- \mathscr{T} v\|_{\widetilde{L}^{2} (I;E^{ct}_{p,1})} \le  \frac{1}{2}  \|  u-  v\|_{\widetilde{L}^{2} (I;E^{ct}_{p,1})}, \quad  u,v \in  \mathscr{D}.
\end{align}

Next, we show that
\begin{align} \label{4.28}
 \sum_{k}\langle k\rangle^{-1} \| 2^{ct|k|}   \Box_k  u\|_{L^{\infty}_{t\in [0, t_0]}L^p_{x}} \lesssim    \|u_0\|_{M^{-1}_{p,1}} + \delta^2.
\end{align}
In order to show \eqref{4.28}, we need the following
\begin{prop} \label{prop5.3a}
Let $1\le p\le \infty$. There exists a constant $c>0$ ($0<c\le 2^{-10}$) such that for $0< t_0 < \infty$ and $I=[0, t_0]$,
\begin{align}
 \sum_{k}\langle k\rangle^{-1} \| 2^{ct|k|}   \Box_k U_2(t) u_0\|_{L^{\infty}_{t\in [0, t_0]}L^p_{x}} & \lesssim    \|u_0\|_{M^{-1}_{p,1}}, \label{4.29} \\
 \sum_{k\neq 0}\langle k\rangle^{-1} \| 2^{ct|k|}   \Box_k  \mathscr{A}_2   \nabla  f\|_{L^{\infty}_{t\in [0, t_0]}L^p_{x}}  & \lesssim \|f\|_{\widetilde{L}^{ 1}(I;E^{ct}_{p, 1})} \label{4.30}
\end{align}
  holds for all $u_0\in M^{-1}_{p,1}$ and  $f\in  \widetilde{L}^{ 1}(I;E^{ct}_{p, 1})$.
  \end{prop}

{\it Proof.} In view of \eqref{4.1}, we immediately have \eqref{4.29} and  \eqref{4.30}.  $\hfill\Box$

Applying Proposition \ref{prop5.3a} and noticing that  $\widetilde{\Box}_0\mathbb{P}{\rm div}: L^p\to L^p$ for all $p\in [1,\infty]$, we have
\begin{align} \label{4.31}
 \sum_{k}\langle k\rangle^{-1} \| 2^{ct|k|}   \Box_k  u\|_{L^{\infty}_{t\in [0, t_0]}L^p_{x}} & \lesssim    \|u_0\|_{M^{-1}_{p,1}} + \|u\otimes u\|_{\widetilde{L}^{ 1}(I;E^{ct}_{p, 1})} \nonumber\\
 & \lesssim    \|u_0\|_{M^{-1}_{p,1}} + 2^{4ct_0} \|u\|^2_{\widetilde{L}^{2}(I;E^{ct}_{p, 1})},
\end{align}
which implies \eqref{4.28}.

We now extend the solution from $I=[0,t_0]$ to $I_1=[t_0,t_1]$ for some $t_1>t_0$ and consider the mapping
\begin{align} \label{4.32}
\mathscr{T}_1: u(t) \to  U_2(t-t_0) u (t_0) + \int^t_{t_0} U_2(t-\tau) \mathbb{P}{\rm div} (u\otimes u)(\tau) d\tau.
\end{align}
\begin{align} \label{4.33}
 \mathscr{D}_1 =\{u\in  \widetilde{L}^{2}(I_1;E^{ct}_{p,1})  : \  \|u\|_{\widetilde{L}^{2} (I_1;E^{ct}_{p,1})} \le \delta_1\},
\end{align}
where $\delta_1$ will be chosen below.  Taking $t=t_0$ in \eqref{4.28}, one has that
 \begin{align} \label{4.34}
 \sum_{k}\langle k\rangle^{-1} \| 2^{ct_0|k|}   \Box_k  u\|_{p} \lesssim    \|u_0\|_{M^{-1}_{p,1}} + \delta^2.
\end{align}
For any $\delta_1>0$,   in view of \eqref{4.34}, we can choose a sufficiently large $J$ such that
\begin{align} \label{4.35}
 C \sum_{|k|>J}\langle k\rangle^{-1} \| 2^{ct_0|k|}   \Box_k  u\|_{p} \le      \delta_1/4.
\end{align}
Hence, in view of Proposition \ref{prop5.2},
\begin{align} \label{4.36}
\|U_2(t-t_0) u (t_0) \|_{\widetilde{L}^{2}(I_1, \{k:|k|>J\};E^{ct}_{p, 1})}   \le C  \sum_{|k|>J} \langle k\rangle^{-1} \| 2^{ct_0|k|}   \Box_k  u\|_{p} \le      \delta_1/4,
\end{align}
and one can choose $t_1>t_0$ verifying $C |I_1|^{1/2} \langle J \rangle (\|u_0\|_{M^{-1}_{p,1}} + \delta^2) \le \delta_1/4$, so,
\begin{align}
\|U_2(t-t_0) u (t_0) \|_{\widetilde{L}^{2}(I_1, \{k:|k| \le J\};E^{ct}_{p, 1})}  &  \le C |I_1|^{1/2} \sum_{|k| \le J}  \| 2^{ct_0|k|}   \Box_k  u\|_{p} \nonumber\\
 &  \le C |I_1|^{1/2} \langle J \rangle  \sum_{|k| \le J} \langle k\rangle^{-1} \| 2^{ct_0|k|}   \Box_k  u\|_{p} \nonumber\\
  &  \le C |I_1|^{1/2} \langle J \rangle (\|u_0\|_{M^{-1}_{p,1}} + \delta^2) \le \delta_1/4.  \label{4.37}
\end{align}
In view of \eqref{4.36} and \eqref{4.37},
\begin{align} \label{4.38}
\|U_2(t-t_0) u (t_0) \|_{\widetilde{L}^{2}(I_1;E^{ct}_{p, 1})}     \le      \delta_1/2,
\end{align}
It follows from Propositions \ref{prop5.2} and \ref{prop5.3}
\begin{align} \label{4.39}
\left\|\int^t_{t_0} U_2(t-\tau) \mathbb{P}{\rm div} (u\otimes u)(\tau) d\tau \right\|_{\widetilde{L}^{2}(I_1;E^{ct}_{p,1})}  &  \le C(1+|I_1|^{1/2}) \|u\otimes u\|_{\widetilde{L}^{1}(I_1;E^{ct}_{p,1})} \nonumber\\
 &  \le   C(1+|I_1|^{1/2}) 2^{4ct_1} \|u\|^2_{\widetilde{L}^{2}(I_1;E^{ct}_{p,1})}.
\end{align}
We can assume that $t_1\le 2$. Hence, collecting \eqref{4.38} and \eqref{4.39}, we have
\begin{align} \label{4.40}
\|\mathscr{T}_1 u\|_{\widetilde{L}^{2}(I_1;E^{ct}_{p,1})}  &  \le \delta_1 /2 +   C  2^{8c} \|u\|^2_{\widetilde{L}^{2}(I_1;E^{ct}_{p,1})}.
\end{align}
Now we can choose $\delta_1>0$ satisfying
\begin{align} \label{4.41}
  C  2^{8c} \delta_1 \le 1/4.
\end{align}
It follows that  $\mathscr{T}_1 u \in \mathscr{D}_1$ for any $ u \in \mathscr{D}_1$. So, we have extended the solution from $[0,t_0] $ to $[t_0,t_1]$.  Noticing that
\begin{align}
 \sum_{k}\langle k\rangle^{-1} \| 2^{ct|k|}   \Box_k U_2(t) u_0\|_{L^{\infty}_{t\in [t_0, t_1]}L^p_{x}} & \lesssim  \sum_{k}\langle k\rangle^{-1} \| 2^{ct_0|k|}   \Box_k  u_0\|_{p}, \label{4.42}
\end{align}
we easily see that the solution can be extended to $[t_1, t_2]$, $[t_2, t_{3}]$,...  and finally find a $T_{\max}>0$ verifying the conclusions of Theorem \ref{NSTh1b}.

\section{Well-posedness of GNS: $\alpha=1/2$ } \label{sect3}

In this section we give some time-space estimates for the solutions of the linear evolution equation
\begin{align}
u_t+(-\Delta)^{1/2}u=f, \ \   u(0,x)= u_0. \label{lin7.1}
\end{align}
 In the following,  we will set up some
estimates for $U_1(t)$
using exponetial decay property of $e^{-t|\xi|}$ together with the dyadic decomposition.
\begin{lem}\label{lem1}
Let $1\leq p \leq \infty$,  $1\leq q \leq \infty$,  $1\leq \gamma \leq
\infty$. Then
\begin{align*}
\|U_1(t)f\|_{\widetilde{L}^{\gamma}(\mathbb{R}^{+};
\dot{B}^{s}_{p, q})}\lesssim \|f\|_{\dot{B}^{s-1/{\gamma}}_{p, q}}.
\end{align*}
\end{lem}
{\it Proof.} By Nikol'skij's inequality,  we easily see that for some $L>n/2$
\begin{align}
\|\Delta_j U_1(t)f\|_{L^p_x}  \leq \|\varphi(\xi)e^{-t2^{j}|\xi|}\|_{M_p}\|f\|_{L^p_x}\lesssim
\|\varphi e^{-t2^{j}|\xi|}\|_{H^{L}}\|f\|_{L^p_x}
 \lesssim
e^{-ct2^{j}}\|f\|_{L^p_x},  \no\label{}
\end{align}
Note that
$\Delta_j=\Delta_j(\Delta_{j-1}+\Delta_{j}+\Delta_{j+1})$,  hence we
have
\begin{align}
\|\Delta_{j}U_1(t)f\|_{L^p_x}\lesssim e^{-ct2^{j}}
\|\Delta_{j}f\|_{L^p_x}.\label{freloc}
\end{align}
Taking $L_t^{\gamma}$ norm on inequality \eqref{freloc},  we get
\begin{align}
\|\Delta_{j}U_1(t)f\|_{L_t^{\gamma}L^p_x}\lesssim
\|e^{-ct2^{j}}\|_{L_t^{\gamma}} \|\Delta_{j}f\|_{L^p_x}\lesssim
2^{-j/{\gamma}}\|\Delta_{j}f\|_{L^p_x}.\label{freloc1}
\end{align}
Multiplying \eqref{freloc1} by $2^{js}$ and taking $l^{q}$ norm,  we
obtain the results as desired.$\hfill\Box$

\begin{lem}\label{lem2}
Let $1\leq p ,q \leq \infty$,  $1\leq \gamma_1 \le \gamma \leq \infty$. Then
\begin{align*}
\|\mathscr{A}_1 f\|_{\widetilde{L}^{\gamma}(\mathbb{R}^{+};
\dot{B}^{s}_{p, q})}\lesssim
\|f\|_{\widetilde{L}^{\gamma_1}(\mathbb{R}^{+};\dot{B}^{s-(1/\gamma+ 1- 1/\gamma_1)}_{p, q})}.
\end{align*}
\end{lem}
{\it Proof.} By Lemma \ref{lem1},
\begin{align} \|\Delta_j\mathscr{A}_1 f\|_{L^p_x} & \leq
\int^t_0 \|\Delta_j U_1(t-\tau)f(\tau, x)\|_{L^p_x}d\tau\no\label{}\\
&\lesssim \int_0^t e^{-c(t-\tau)2^{j}}\|\Delta_j
f(\tau)\|_{L^p_x}d\tau,  \no\label{}
\end{align}
Taking
$L^{\gamma}$ norm on time variable and using Young inequality,
\begin{align} & \|\Delta_j\mathscr{A}_1 f\|_{L_t^{\gamma}L^p_x}
 \lesssim 2^{-j(1+1/\gamma-1/\gamma_1)}\|\Delta_j f\|_{L_t^{\gamma_1}L^p_x}.
\label{freloc2}
\end{align}
Multiplying \eqref{freloc2} by $2^{js}$ and taking $l^{q}$ norm,  we
obtain the results.$\hfill\Box$ \\

Below, we prove the existence part of Theorem \ref{Th1} and the analyticity part is left into the next section.   We  divide the proof into two cases.

 Case I: $1\leq p <\infty$.
We construct the resolution space as follows:
\begin{align*}
\mathscr{D}=\{u\in \mathscr{S}'(\mathbb{R}^{+}\times
\mathbb{R}^{n}): \|u\|_{\widetilde{L}^{\infty}(\mathbb{R}^{+};
\dot{B}_{p, 1}^{n/p})}\leq 2C\|u_0\|_{\dot{B}_{p, 1}^{n/p}}\}
\end{align*}
with metric
$$
d(u, v)=\|u-v\|_{\widetilde{L}^{\infty}(\mathbb{R}^{+}; \dot{B}_{p,
1}^{n/p})}.
$$
Assume that $u_0\in \dot{B}^{n/p}_{p, 1}$ is sufficiently small such that
$4C\|u_0\|_{\dot{B}^{n/p}_{p, 1}}=1/2$,  where $C$ is the largest
constant appeared in the following inequalities. Considering the mapping
$$
\mathscr{T}: u(t)\rightarrow
U_1(t)u_0+\mathscr{A}_1\mathbb{P}\mathrm{div}(u\otimes u),
$$
we will verify that $\mathscr{T}:(\mathscr{D},  d)\rightarrow
(\mathscr{D},  d)$ is a contraction mapping. In fact,  by Lemma
\ref{lem1} and Lemma \ref{lem2},  we have
\begin{align}
\|\mathscr{T} u\|_{\widetilde{L}^{\infty}(\mathbb{R}_{+};\dot{B}_{p,
1}^{n/p})}
& \lesssim \|u_0\|_{\dot{B}_{p, 1}^{n/p}}+\|u\otimes
u\|_{\widetilde{L}^{\infty}(\mathbb{R}_{+};\dot{B}_{p, 1}^{n/p})}.
\label{contra}
\end{align}
We will use the decomposition \eqref{para1} to deal with the
estimates of $\|u\otimes
u\|_{\widetilde{L}^{\infty}(\mathbb{R}_{+};\dot{B}_{p, 1}^{n/p})}$.
We have
\begin{align}
\|\Delta_j(u\otimes u)\|_{L_x^p} &=\|\sum_{j\leq
k+3}\Delta_j(S_{k-1}u\Delta_ku+S_k u \Delta_k u)\|_{L_x^p}
\nonumber\\
&\lesssim \sum_{k\geq j-3}\|S_k
u\|_{L_x^{\infty}}\|\Delta_k u\|_{L_x^p}
\no\label{}\\
&\lesssim \sum_{k\geq j-3}\sum_{l\leq k}\|\Delta_{l}
u\|_{L_x^{\infty}}\|\Delta_k u\|_{L_x^p} \nonumber \\
& \lesssim \sum_{k\geq j-3}\sum_{l\leq k}2^{ln/p}\|\Delta_{l}
u\|_{L_x^{p}}\|\Delta_k u\|_{L_x^p},  \label{contraa}
\end{align}
where we used Bernstein's inequality in the last estimate. Hence
taking $L_t^{\infty}$ norm to \eqref{contraa},  we obtain
\begin{align}
\|\Delta_j(u\otimes u)\|_{L_t^{\infty}L_x^p} & \lesssim \sum_{k\geq
j-3}\sum_{l\leq k}2^{ln/p}\|\Delta_{l}
u\|_{L_t^{\infty}L_x^{p}}\|\Delta_k
u\|_{L_t^{\infty}L_x^p}\nonumber\\
&\lesssim \|u\|_{\widetilde{L}^{\infty}(\mathbb{R}_{+};\dot{B}_{p,
1}^{n/p})} \sum_{k\geq j-3}\|\Delta_k u\|_{L_t^{\infty}L_x^p}.
\label{contra1}
\end{align}
Multiplying \eqref{contra1} by $2^{jn/p}$,
\begin{align}
2^{jn/p}\|\Delta_j(u\otimes u)\|_{L_t^{\infty}L_x^p}& \lesssim
\|u\|_{\widetilde{L}^{\infty}(\mathbb{R}_{+};\dot{B}_{p,
1}^{n/p})} \sum_{k\geq j-3}2^{-(k-j)n/p}2^{kn/p}\|\Delta_k
u\|_{L_t^{\infty}L_x^p}. \label{contra2}
\end{align}
Taking $l^1$ norm in \eqref{contra2} and then using Young's inequality,
\begin{align}
\|u\otimes u\|_{\widetilde{L}^{\infty}(\mathbb{R}_{+};\dot{B}_{p,
1}^{n/p})}
& \lesssim \|u\|^2_{\widetilde{L}^{\infty}(\mathbb{R}_{+};\dot{B}_{p,
1}^{n/p})} .\label{contra3a}
\end{align}
Inserting \eqref{contra3a} into \eqref{contra},  we get
\begin{align}
\|\mathscr{T} u\|_{\widetilde{L}^{\infty}(\mathbb{R}_{+};\dot{B}_{p,
1}^{n/p})} & \lesssim \|u_0\|_{\dot{B}_{p,
1}^{n/p}}+\|u\|^2_{\widetilde{L}^{\infty}(\mathbb{R}_{+};\dot{B}_{p,
1}^{n/p})}
\nonumber\\
&\leq C \|u_0\|_{\dot{B}_{p, 1}^{n/p}} + 4C^2
\|u_0\|^2_{\dot{B}^{n/p}_{p, 1}} <2C\|u_0\|_{\dot{B}_{p, 1}^{n/p}},
\nonumber
\end{align}
which means that $\mathscr{T}: \mathscr{D}\rightarrow \mathscr{D}$. Similarly,
\begin{align}
\|\mathscr{T} u-\mathscr{T}
v\|_{\widetilde{L}^{\infty}(\mathbb{R}_{+};\dot{B}_{p, 1}^{n/p})}
& <
\frac{1}{2}\|u-v\|_{\widetilde{L}^{\infty}(\mathbb{R}_{+};\dot{B}_{p,
1}^{n/p})}. \nonumber
\end{align}
Hence there exists a $u\in
\widetilde{L}^{\infty}(\mathbb{R}_{+};\dot{B}_{p, 1}^{n/p})$
satisfying
$$
u(t)=U_1(t)u_0+\mathscr{A}_1\mathbb{P}\mathrm{div}(u\otimes u)
$$
with
$$
\|u\|_{\widetilde{L}^{\infty}(\mathbb{R}_{+};\dot{B}_{p,
1}^{n/p})}\lesssim \|u_0\|_{\dot{B}_{p, 1}^{n/p}}.
$$
In particular,  we have
$$
\|u\|_{L^{\infty}(\mathbb{R}_{+};\dot{B}_{p, 1}^{n/p})}\lesssim
\|u\|_{\widetilde{L}^{\infty}(\mathbb{R}_{+};\dot{B}_{p,
1}^{n/p})}\lesssim \|u_0\|_{\dot{B}_{p, 1}^{n/p}}.
$$
The uniqueness is similar. \\

Case II: $p=\infty$.  We now consider the following resolution
space:
\begin{align*}
\mathscr{D}_1=\{u\in \mathscr{S}'(\mathbb{R}^{+}\times
\mathbb{R}^{n}): \|u\|_{\widetilde{L}^{\infty}(\mathbb{R}^{+};
\dot{B}^{0}_{\infty,1})\cap \widetilde{L}^{1}(\mathbb{R}^{+};
\dot{B}^{1}_{\infty,1})}\leq 2C\|u_0\|_{\dot{B}^{0}_{\infty,1}}\}
\end{align*}
with metric
$$
d_1(u,v)=\|u-v\|_{\widetilde{L}^{\infty}(\mathbb{R}^{+};
\dot{B}^{0}_{\infty,1})\cap \widetilde{L}^{1}(\mathbb{R}^{+};
\dot{B}^{1}_{\infty,1})}.
$$
Suppose that $u_0\in \dot{B}^{0}_{\infty,1}$ is sufficiently small
such that $4C\|u_0\|_{\dot{B}^{0}_{\infty,1}}=1/2$. We will verify that $\mathscr{T}:(\mathscr{D}_1, d_1)\rightarrow
(\mathscr{D}_1, d_1)$ is a contractive map. Using Lemma \ref{lem1} and
Lemma \ref{lem2},
\begin{align}
\|\mathscr{T} u\|_{\widetilde{L}^{\infty}(\mathbb{R}^{+};
\dot{B}^{0}_{\infty,1})\cap \widetilde{L}^{1}(\mathbb{R}^{+};
\dot{B}^{1}_{\infty,1})}
&\lesssim \|u_0\|_{\dot{B}^{0}_{\infty,1}}+\|u\otimes
u\|_{\widetilde{L}^{1}(\mathbb{R}_{+};\dot{B}_{\infty,1}^{1})},
\label{contra3}
\end{align}
From \eqref{contraa}, we see that
\begin{align} &
\|\Delta_j(u\otimes u)\|_{L_x^{\infty}} \lesssim \sum_{k\geq
j-3}\sum_{l\leq k}\|\Delta_{l} u\|_{L_x^{\infty}}\|\Delta_k
u\|_{L_x^{\infty}}, \label{contra4}
\end{align}
Taking $L_t^{1}$ norm on \eqref{contra4}, we obtain
\begin{align}
\|\Delta_j(u\otimes u)\|_{L_t^{1}L_x^{\infty}}
 \lesssim
\|u\|_{\widetilde{L}^{\infty}(\mathbb{R}_{+};\dot{B}_{\infty,1}^{0})}
\sum_{k\geq j-3}\|\Delta_k u\|_{L_t^{1}L_x^{\infty}}.
\label{contra5}
\end{align}
Multiplying \eqref{contra5} by $2^{j}$ and then
taking $l^{1}$ norm on \eqref{contra5},
\begin{align}
\|u\otimes
u\|_{\widetilde{L}^{1}(\mathbb{R}_{+};\dot{B}^{1}_{\infty,1})}
\lesssim
\|u\|_{\widetilde{L}^{\infty}(\mathbb{R}_{+};\dot{B}_{\infty,1}^{0})}\left\|\sum_{k\geq
j-3}2^{-(k-j)}2^{k}\|\Delta_k
u\|_{L_t^{1}L_x^{\infty}}\right\|_{l^{1}}. \nonumber
\end{align}
Therefore, in view of Young's inequality, we have
\begin{align}&\|u\otimes
u\|_{\widetilde{L}^{1}(\mathbb{R}_{+};\dot{B}^{1}_{\infty,1})}
\lesssim
\|u\|_{\widetilde{L}^{\infty}(\mathbb{R}_{+};\dot{B}_{\infty,1}^{0})}
\|u\|_{\widetilde{L}^{1}(\mathbb{R}_{+};\dot{B}^{1}_{\infty,1})}.
\label{contra7}
\end{align}
Combining \eqref{contra7} with \eqref{contra3}, we get for any $u\in \mathscr{D}_1$,
\begin{align} &
\|\mathscr{T} u\|_{\widetilde{L}^{\infty}(\mathbb{R}^{+};
\dot{B}^{0}_{\infty,1})\cap \widetilde{L}^{1}(\mathbb{R}^{+};
\dot{B}^{1}_{\infty,1})}\lesssim
\|u_0\|_{\dot{B}^{0}_{\infty,1}}+\|u\|_{\widetilde{L}^{\infty}(\mathbb{R}_{+};\dot{B}_{\infty,1}^{0})}
\|u\|_{\widetilde{L}^{1}(\mathbb{R}_{+};\dot{B}^{1}_{\infty,1})}\nonumber\\
&\leq C \|u_0\|_{\dot{B}^{0}_{\infty,1}} + 4C^2
\|u_0\|_{\dot{B}^{0}_{\infty,1}}\|u_0\|_{\dot{B}^{0}_{\infty,1}}<2C\|u_0\|_{\dot{B}^{0}_{\infty,1}},\nonumber
\end{align}
which means that $\mathscr{T}: \mathscr{D}_1\rightarrow
\mathscr{D}_1$ and there exists a solution $u\in
\widetilde{L}^{\infty}(\mathbb{R}^{+}; \dot{B}^{0}_{\infty,1})\cap
\widetilde{L}^{1}(\mathbb{R}^{+}; \dot{B}^{1}_{\infty,1})$ to
\eqref{1.1} satisfying
$$
\|u\|_{\widetilde{L}^{\infty}(\mathbb{R}^{+};
\dot{B}^{0}_{\infty,1})\cap \widetilde{L}^{1}(\mathbb{R}^{+};
\dot{B}^{1}_{\infty,1})}\lesssim \|u_0\|_{\dot{B}^{0}_{\infty,1}}.
$$
In particular, we have
$$
\|u\|_{L^{\infty}(\mathbb{R}_{+};\dot{B}^{0}_{\infty,1})}\lesssim
\|u\|_{\widetilde{L}^{\infty}(\mathbb{R}_{+};\dot{B}^{0}_{\infty,1})}\lesssim
\|u_0\|_{\dot{B}^{0}_{\infty,1}}.
$$
The uniqueness is similar, we omit the details here.$\hfill\Box$

 \begin{rem} \rm Since $\widetilde{L}^{\infty}(\mathbb{R}^{+};
\dot{B}^{0}_{\infty,1})$ is not a Banach algebra,  the method for the nonlinear estimates in
case $1\leq p<\infty$ doesn't work for the case
 $p=\infty$.
 \end{rem}

\section{Analyticity of GNS: $\alpha=1/2$} \label{sectGNS}

 First,  we show that there exists $c>0$ (say $0<c \le 2^{-5}$) such that
\begin{align} \label{7.1}
\|\Box_k    U_1(t) f\|_{p}  \lesssim  e^{-2ct|k|} \|\Box_k f \|_{p}
\end{align}
holds for all $f\in L^p$ and $k\in \mathbb{Z}^n$. Let $\tilde{\sigma}:\mathbb{R}^n \to [0,1]$ be  the same cut-off function as in Section \ref{sectNS}. In view of Nikol'skij's inequality, we have
\begin{align}
\|\Box_k    U_1(t) f\|_{p}
& \lesssim    \|  \tilde{\sigma} (\xi-k) e^{- t|\xi|}\|_{M_p} \|\Box_k f \|_{p} \nonumber\\
& \lesssim \| \tilde{\sigma}  e^{- t|\xi+k|}\|_{H^L} \|\Box_k f \|_{p}
\end{align}
for any $ L>n(1/(1\wedge p)-1/2).$ Similar to \eqref{4.4} and \eqref{4.5}, we have
\begin{align}
 | \partial^L_{\xi_i} (\tilde{\sigma}  e^{- t|\xi+k|}) |   \lesssim  e^{- t|\xi+k|/2} \chi_{[-7/8,7/8]^n}(\xi),
\end{align}
we immediately have for $|k|  \ge 1$,
\begin{align} \label{7.6}
  \| \tilde{\sigma} (\xi-k) e^{- t|\xi|}\|_{M_p}   \lesssim  e^{- t|k|/16}.
\end{align}
We easily see that  \eqref{7.6} also holds for $k=0$. Hence,   we have shown that

\begin{prop} \label{prop6.1}
Let $0<p\le \infty$. There exists a constant $c>0$ (say $0< c \le 2^{-5}$) such that
\begin{align} \label{7.7}
\|\Box_k    U_1(t) f\|_{p}  \lesssim  2^{-2ct|k|} \|\Box_k f \|_{p}
\end{align}
uniformly holds for all $k\in \mathbb{Z}^n$ and $f\in L^p$.
  \end{prop}

Now let $0<t_0 \le \infty$.  We consider the estimate of $  \mathscr{A}_1  f $. Applying Proposition \ref{prop6.1},
we have
\begin{align} \label{7.8}
 \|\Box_k  \mathscr{A}_1  f\|_{L^p_{x}} \lesssim \int^t_0   2^{-2c(t-\tau)|k|} \|\Box_k f (\tau)\|_{L^p_{x}} d\tau.
\end{align}
It follows from \eqref{7.8} that
\begin{align} \label{7.9}
 2^{ct|k|} \|\Box_k  \mathscr{A}_1  f\|_{L^p_{x}} \lesssim \int^t_0  2^{- c(t-\tau)|k|} 2^{ c\tau |k|}\|\Box_k f (\tau)\|_{L^p_{x}} d\tau.
\end{align}
For $|k| \ge 1$ and $\tilde{q} \ge \tilde{q}_1$,  using Young's inequality, we have from \eqref{7.9} that
\begin{align} \label{7.10}
 \| 2^{ct|k|} \Box_k  \mathscr{A}_1  f\|_{L^{\tilde{q}}_{ t\in [0, t_0)}L^p_{x}} \lesssim  \langle k\rangle ^{1/\tilde{q}_1-1/\tilde{q}- 1} \|2^{ c\tau |k|}\Box_k f (\tau)\|_{L^{\tilde{q}_1}_{ \tau\in [0, t_0)}L^p_{x}}.
\end{align}
So,  taking the sequence $l^q$ norm in both sides of \eqref{7.10},  we have
\begin{prop} \label{prop6.2}
Let $1\le p\le \infty$, $1\le \tilde{q}_1 \le \tilde{q} \le \infty$. There exists a constant $c>0$ (say, $c \le 2^{-5}$) such that for $0< t_0 \le \infty$, $I=[0, t_0)$ and $\mathbb{A} \subset  \mathbb{Z}^n \setminus \{0\}$,
\begin{align} \label{7.11}
\|\mathscr{A}_1 f\|_{\widetilde{L}^{\tilde{q}}(I, \mathbb{A};E^{ct}_{p, q})}   \lesssim  \sum_{k\in \mathbb{A}} \langle k\rangle ^{1/\tilde{q}_1-1/\tilde{q}- 1} \|2^{ ct |k|}\Box_k f \|_{L^{\tilde{q}_1}_{ t\in [0, t_0)}L^p_{x}}
\end{align}
  holds for all   $f\in  \widetilde{L}^{{{\tilde{q}}}}(I;E^{ct}_{p, q})$. In particular, we have
\begin{align} \label{7.12}
\|\nabla\mathscr{A}_1 f\|_{\widetilde{L}^{1}(I, \mathbb{A};E^{ct}_{\infty, 1})} &  \lesssim  \| f\|_{\widetilde{L}^{1}(I, \mathbb{A};E^{ct}_{\infty, 1})}, \\
\|\mathscr{A}_1 f\|_{\widetilde{L}^{\infty}(I, \mathbb{A};E^{ct}_{\infty, 1})} &  \lesssim  \| f\|_{\widetilde{L}^{1}(I, \mathbb{A};E^{ct}_{\infty, 1})}. \label{7.13}
\end{align}
\end{prop}

{\it Proof of iv) of Theorem \ref{Th1}.} Denote $I=[0,1]$. Put
\begin{align} \label{7.14}
\|u\|_{X} =  \|u\|_{\widetilde{L}^{1} (I ; \dot{B}^{1}_{\infty,1}) \cap \widetilde{L}^{\infty} (I ; \dot{B}^{0}_{\infty,1})  }, \quad \|u\|_Y= \|\nabla u\|_{\widetilde{L}^{1} (I ; E^{ct}_{\infty,1}) }+ \|u\|_{\widetilde{L}^{\infty} (I ; E^{ct}_{\infty,1})  }
\end{align}
and
\begin{align} \label{7.15}
 \mathscr{D}  =\{u   : \  \|u\|_{X} + \ \| u\|_{Y}  \le \delta\}, \quad d(u,v)=\|u-v\|_{X\cap Y}.
\end{align}
We will show that
\begin{align} \label{7.16}
\mathscr{T}: u(t) \to  U_1(t) u_0 + \mathscr{A}_1 \mathbb{P}{\rm div} (u\otimes u)
\end{align}
is a contraction mapping from $(\mathscr{D},d)$ into itself. In Section \ref{sect3}, we have shown that
\begin{align} \label{7.17}
\|\mathscr{T}  u\|_X  \le C \|u_0\|_{\dot{B}^{0}_{\infty,1}} + C \|u\|^2_X.
\end{align}
So, it suffices to estimate $\|\mathscr{T}  u\|_Y$.
   Put $c=2^{-5}$. By Proposition  \ref{prop6.1},
\begin{align} \label{7.18}
\|\nabla U_1(t) u_0\|_{\widetilde{L}^{1}(I;E^{ct}_{\infty, 1})} + \|U_1(t) u_0\|_{\widetilde{L}^{\infty}(I;E^{ct}_{\infty, 1})}   \le C \|u_0\|_{M^0_{\infty,1}}
\end{align}
Since $\|\Box_0 f\|_\infty \lesssim  \|f\|_\infty \lesssim \|f\|_{\dot{B}^{0}_{\infty,1}} $, we have from Section \ref{sect3} that
\begin{align} \label{7.19}
 & \|\Box_0 \nabla \mathscr{A}_1 \mathbb{P}{\rm div} (u\otimes u) \|_{ {L}^{1}_{t\in I}L^\infty_x} +\|\Box_0 \mathscr{A}_1 \mathbb{P}{\rm div} (u\otimes u) \|_{ {L}^{\infty}_{t\in I}L^\infty_x} \nonumber\\
 & \lesssim \|   \mathscr{A}_1 \mathbb{P}{\rm div} (u\otimes u) \|_{\widetilde{L}^{1}(I; \dot{B}^{1}_{\infty,1}) \cap \widetilde{L}^{\infty}(I; \dot{B}^{0}_{\infty,1}) }   \lesssim   \|u\|^2_X.
\end{align}
Therefore, it suffices to consider the estimate of $\nabla \mathscr{A}_1 \mathbb{P}{\rm div} (u\otimes u)$  in the space $\widetilde{L}^{1}(I, \mathbb{Z}^n_{\star}; E^{ct}_{\infty, 1})$ and $ \mathscr{A}_1 \mathbb{P}{\rm div} (u\otimes u)$  in the space $\widetilde{L}^{\infty}(I, \mathbb{Z}^n_{\star}; E^{ct}_{\infty, 1})$ for $\mathbb{Z}^n_{\star}=\mathbb{Z}^n\setminus \{0\}$.
It follows from Propositions \ref{prop6.2} and \ref{prop5.3} that
\begin{align} \label{7.20}
\|\nabla \mathscr{A}_1 \mathbb{P}{\rm div} (u\otimes u) & \|_{\widetilde{L}^{1}(I, \mathbb{Z}^n_{\star};E^{ct}_{\infty, 1})} + \|\mathscr{A}_1 \mathbb{P}{\rm div} (u\otimes u)\|_{\widetilde{L}^{\infty}(I, \mathbb{Z}^n_{\star};E^{ct}_{\infty,1})} \nonumber\\
  &  \lesssim  \|\nabla (u\otimes u)\|_{\widetilde{L}^{1} (I;E^{ct}_{\infty,1})}
    \lesssim \|\nabla  u\|_{\widetilde{L}^{1} (I;E^{ct}_{\infty,1})} \|u\|_{\widetilde{L}^{\infty} (I;E^{ct}_{\infty,1})}.
\end{align}
Hence, collecting \eqref{7.18}--\eqref{7.20}, we have
\begin{align} \label{7.21}
\|\mathscr{T} u\|_{Y}  &  \le C  \|u_0\|_{M^0_{\infty,1}} + C \|\nabla  u\|_{\widetilde{L}^{1} (I;E^{ct}_{\infty,1})} \|u\|_{\widetilde{L}^{\infty} (I;E^{ct}_{\infty,1})}.
\end{align}
By \eqref{7.17} and \eqref{7.21}
\begin{align} \label{7.21a}
\|\mathscr{T}  u\|_X + \|\mathscr{T} u\|_{Y}  \le C \|u_0\|_{\dot{B}^{0}_{\infty,1} \cap M^0_{\infty,1}}   + C \|u\|^2_X  + C \|u\|^2_Y.
\end{align}
Assume that
 \begin{align} \label{7.22}
C\delta \le 1/4, \quad C \|u_0\|_{\dot{B}^{0}_{\infty,1} \cap M^0_{\infty,1}} \le \delta/2.
\end{align}
Using the above contraction mapping argument, we obtain that there exists a solution of the  integral equation $\mathscr{T}u = u$ in the space $X\cap Y$. In particular, we have
\begin{align} \label{7.23}
\|u\|_{\widetilde{L}^{\infty}(0,1;E^{ct}_{\infty,1})} \le \delta.
\end{align}
In particular, we have for $u(1)=u(1,x)$,
\begin{align} \label{7.24}
\|u(1)\|_{ E^{c}_{\infty,1} } \le \delta, \quad c=2^{-5}.
\end{align}
We now extend the regularity of solutions. Let us consider the mapping
\begin{align} \label{7.25}
\mathscr{T}_1: u(t) \to  U_1(t-1) u (1) + \int^t_1 U_1(t-\tau) \mathbb{P}{\rm div} (u\otimes u)(\tau) d\tau.
\end{align}
Let $I_1=[1, \infty)$ and
\begin{align} \label{7.26}
\|u\|_{X_1} =  \|u\|_{\widetilde{L}^{1} (I_1 ; \dot{B}^{1}_{\infty,1}) \cap \widetilde{L}^{\infty} (I_1 ; \dot{B}^{0}_{\infty,1})  }, \ \|u\|_{Y_1}= \|\nabla u\|_{\widetilde{L}^{1} (I_1 ; E^{c}_{\infty,1}) }+ \|u\|_{\widetilde{L}^{\infty} (I_1 ; E^{c}_{\infty,1})},
\end{align}
\begin{align} \label{7.27}
 \mathscr{D}_1  =\{u   : \  \|u\|_{X_1} + \ \| u\|_{Y_1}  \le \delta\}, \quad d_1(u,v)=\|u-v\|_{X_1\cap Y_1}.
\end{align}
We use the same notation $\mathbb{Z}^n_{\star} = \mathbb{Z}^n \setminus \{0\}$ as above.  Using the same way as in \eqref{7.19}, we see that
\begin{align} \label{7.28}
  \|u\|_{Y_1} \lesssim \|u\|_{X_1} + \|\nabla u\|_{\widetilde{L}^{1} (I_1, \mathbb{Z}^n_{\star} ; E^{c}_{\infty,1}) }+ \|u\|_{\widetilde{L}^{\infty} (I_1, \mathbb{Z}^n_{\star} ; E^{c}_{\infty,1})}.
\end{align}
In section \ref{sect3}, we have shown that
\begin{align} \label{7.29}
\|\mathscr{T}_1  u\|_{X_1}  \le C \|u_0\|_{\dot{B}^{0}_{\infty,1}} + C \|u\|^2_{X_1}.
\end{align}
So, it suffices to estimate $ \nabla\mathscr{T}_1  u $ in the space $\widetilde{L}^{1} (I_1, \mathbb{Z}^n_{\star} ; E^{c}_{\infty,1})$ and $ \mathscr{T}_1  u $ in the space $\widetilde{L}^{\infty} (I_1, \mathbb{Z}^n_{\star} ; E^{c}_{\infty,1})$.
   Put $c=2^{-5}$. By Proposition  \ref{prop6.1},
\begin{align} \label{7.30}
\|\nabla  U_1(t-1) u (1) \|_{\widetilde{L}^{1}(I_1, \mathbb{Z}^n_{\star} ; E^{c }_{\infty, 1})} + \|U_1(t-1) u (1)\|_{\widetilde{L}^{\infty}(I_1, \mathbb{Z}^n_{\star} ;E^{c }_{\infty, 1})}   \le C \|u(1)\|_{E^c_{\infty,1}}.
\end{align}
By \eqref{7.8} and Young's inequality,
\begin{align} \label{7.31}
\left\|\nabla \int^t_1 U_1(t-\tau) f(\tau) d\tau \right\|_{\widetilde{L}^{1}(I_1, \mathbb{Z}^n_{\star}; E^{c}_{\infty,1})}  \lesssim  \| f\|_{\widetilde{L}^{1}(I_1, \mathbb{Z}^n_{\star}; E^{c}_{\infty,1})}, \\
 \left\| \int^t_1 U_1(t-\tau) f(\tau) d\tau \right\|_{\widetilde{L}^{\infty}(I_1, \mathbb{Z}^n_{\star}; E^{c}_{\infty,1})} \lesssim  \| f\|_{\widetilde{L}^{1}(I_1, \mathbb{Z}^n_{\star}; E^{c}_{\infty,1})}. \label{7.32}
\end{align}
Therefore, by \eqref{7.31} and Proposition \ref{prop5.3} we have
\begin{align} \label{7.33}
\left\|\nabla \int^t_1 U_1(t-\tau) \mathbb{P}{\rm div} (u\otimes u)(\tau) d\tau \right\|_{\widetilde{L}^{1}(I_1, \mathbb{Z}^n_{\star}; E^{c}_{\infty,1})}
& \lesssim  \| \mathbb{P}{\rm div} (u\otimes u)\|_{\widetilde{L}^{1}(I_1, \mathbb{Z}^n_{\star}; E^{c}_{\infty,1})} \nonumber \\
&  \lesssim  \| {\rm div} (u\otimes u)\|_{\widetilde{L}^{1}(I_1; E^{c}_{\infty,1})} \nonumber \\
&  \lesssim  \|\nabla u\|_{\widetilde{L}^{1}(I_1; E^{c}_{\infty,1})}  \|u\|_{\widetilde{L}^{\infty}(I_1; E^{c}_{\infty,1})} \nonumber \\
&  \lesssim  \|u\|^2_{Y_1}.
\end{align}
By \eqref{7.32} we see that $\left\|\int^t_1 U_1(t-\tau) \mathbb{P}{\rm div} (u\otimes u)(\tau) d\tau \right\|_{\widetilde{L}^{\infty}(I_1, \mathbb{Z}^n_{\star}; E^{c}_{\infty,1})}$ also has the up-bound as in \eqref{7.33}. Hence, collecting \eqref{7.24}, \eqref{7.29}, \eqref{7.30} and \eqref{7.33}, we have
\begin{align} \label{7.34}
  \|\mathscr{T}_1 u\|_{Y_1} & \lesssim  C \|u_0\|_{\dot{B}^{0}_{\infty,1}} + C \|u(1)\|_{ E^{c}_{\infty,1} } + C \|u\|^2_{X_1} + C \|u\|^2_{Y_1} \nonumber\\
& \lesssim  C \delta + C \|u\|^2_{X_1 \cap Y_1}.
\end{align}
Using the contraction mapping argument, we can show that the integral equation $\mathscr{T}_1 u = u$ has a solution in $\mathscr{D}_1$.
The left part of the proof is standard and we omit the details. $\hfill\Box$\\

{\it Proof of iii) of Theorem \ref{Th1}.}
  In view of $e^{ -\frac{3}{4}t|\xi| +  \frac{1}{2n} t|\xi|_1} \in M_p$,  one easily sees that
\begin{align}
\|\Delta_j e^{t\Lambda/2n} U_1(t) u_0\|_{p} \le    \|\Delta_j U_1(t/4)
u_0\|_{p}  \lesssim   e^{-ct2^{j}} \|\Delta_j
u_0\|_{p}.   \label{GNSfreloca}
\end{align}
By \eqref{GNSfreloca},
\begin{align} \label{6.35}
\|\Delta_j e^{t\Lambda/2n} \mathscr{A}_1 f\|_{p}
   \lesssim \int_0^t e^{-c(t-\tau)2^{j}}
 \|\Delta_j  e^{\tau \Lambda/2n} f(\tau)\|_{p} d\tau.
\end{align}
Using Young's inequality,
\begin{align} \label{6.36}
\|\Delta_j e^{t\Lambda/2n} \mathscr{A}_1 f\|_{L^\infty_{t\in I}L^p_x}
    \lesssim   2^{-  j }
 \|\Delta_j  e^{t \Lambda/2n} f \|_{L^{\infty}_{t \in I}L^p_x }.
\end{align}
By \eqref{6.36}, we have
\begin{align}\label{6.37}
\|\nabla \mathscr{A}_1 f\|_{\widetilde{L}^{\infty}(I;
e^{t\Lambda/2n} \dot{B}^{s}_{p,q})} \lesssim
\|f\|_{\widetilde{L}^{\infty}(I;e^{t\Lambda/2n} \dot{B}^{s}_{p,q})}.
\end{align}
We consider the mapping
\begin{align} \label{6.38}
\mathscr{T}: u(t) \to U_1(t) u_0 + \mathscr{A}_1 \mathbb{P}{\rm div} (u\otimes u)
\end{align}
in the metric space
\begin{align} \label{6.39}
\mathscr{D} = \{u: \|u\|_{\widetilde{L}^{\infty}(\mathbb{R}^+; e^{t\Lambda/2n} \dot{B}^{n/p}_{p,1})} \le \delta\}, \quad d(u,v) =\|u-v\|_{\widetilde{L}^{\infty}(\mathbb{R}^+; e^{t\Lambda/2n} \dot{B}^{n/p}_{p,1})}.
\end{align}
By \eqref{GNSfreloca},  \eqref{6.37} and Lemma \ref{NSlem3} (it is easy to see that one can replace $\sqrt{t}$ by $t/2n$ in Lemma \ref{NSlem3}),  for any $u, v\in \mathscr{D}$,
\begin{align} \label{6.40}
   \|\mathscr{T} u\|_{\widetilde{L}^{\infty}(\mathbb{R}^+; e^{t\Lambda/2n} \dot{B}^{n/p}_{p,1})} & \lesssim  \|U_1(t) u_0\|_{\widetilde{L}^{\infty}(\mathbb{R}^+; e^{t\Lambda/2n} \dot{B}^{n/p}_{p,1}) } +  \|u\otimes u\|_{\widetilde{L}^{\infty}(\mathbb{R}^+; e^{t\Lambda/2n} \dot{B}^{n/p}_{p,1})} \nonumber\\
 &  \lesssim  \|u_0\|_{\dot{B}^{n/p}_{p,1}} +  \| u\|^2_{\widetilde{L}^{\infty}(\mathbb{R}^+; e^{t\Lambda/2n} \dot{B}^{n/p}_{p,1})}  \nonumber\\
 &  \lesssim  \|u_0\|_{\dot{B}^{n/p}_{p,1}} +  \delta^2, \\
\|\mathscr{T} u - \mathscr{T}v\|&_{\widetilde{L}^{\infty}(\mathbb{R}^+; e^{t\Lambda/2n } \dot{B}^{n/p}_{p,1})}   \lesssim  \delta \|u-v\|_{\widetilde{L}^{\infty}(\mathbb{R}^+; e^{t\Lambda/2n} \dot{B}^{n/p}_{p,1})}. \label{6.41}
\end{align}
Applying the standard  argument, we can show that $\mathscr{T}$ is a contraction mapping from $(\mathscr{D},d) \to (\mathscr{D},d)$.  So, there is a $u\in \mathscr{D}$ satisfying $\mathscr{T} u=u$.  The left part of the proof is standard and we omit it. $\hfill\Box$\\

\section{Analyticity of GNS: $1/2<\alpha <1$} \label{GNS7}

\begin{lem}\label{7-1}
Let $1/2<\alpha <1$, $1 < p < \infty$,  $1\leq q \leq \infty$,  $1\leq \gamma
\leq \infty$. Then
\begin{align} \label{N7.1}
\|U_{2\alpha}(t) u_0\|_{\widetilde{L}^{\gamma}(\mathbb{R}^{+}; e^{t^{1/2\alpha}\Lambda} \dot{B}^{s}_{p,
q})}\lesssim \|u_0 \|_{\dot{B}^{s-2\alpha/{\gamma}}_{p, q}}.
\end{align}
\end{lem}
{\it Sketch of Proof.} We have
\begin{align}
\|\Delta_j U_{2\alpha}(t)f\|_{L^p_x} \lesssim e^{-ct2^{2\alpha j}}\|\Delta_j
f\|_{L^p_x},  \label{N7.2}
\end{align}
In view of \eqref{N7.2},  one easily sees that
\begin{align}
\|\Delta_j e^{t^{1/2\alpha}\Lambda} U_{2\alpha}(t) u_0\|_{p}  \lesssim e^{-ct2^{2\alpha j}}  \|\Delta_j
u_0\|_{p}.   \label{N7.3}
\end{align}
Taking $L_t^{\gamma}$ norm on inequality \eqref{N7.3},
\begin{align}
\|\Delta_j e^{ t^{1/2\alpha} \Lambda} U_{2\alpha}(t) u_0\|_{L^\gamma_tL^p_x} \lesssim
2^{-2\alpha j/{\gamma}}\|\Delta_{j} u_0\|_{p}.\label{N7.4}
\end{align}
Taking sequence $l^q$ norms in \eqref{N7.4}, we have the result, as desired. $\hfill\Box$

\begin{lem}\label{7-2}
Let $1/2<\alpha <1$, $1<p<\infty$, $1\leq   q\leq \infty$, $1\leq \gamma_1 \leq \gamma \le \infty$, $I= [0,T)$, $T\le \infty$. Then
\begin{align}
\|\mathscr{A}_{2\alpha} f\|_{\widetilde{L}^{\gamma}(I;
e^{t^{1/2\alpha} \Lambda} \dot{B}^{s}_{p,q})}\lesssim
\|f\|_{\widetilde{L}^{\gamma_1}(I;e^{t^{1/2\alpha} \Lambda} \dot{B}^{s-2\alpha(1+1/\gamma -1/\gamma_1)}_{p,q})}.
\end{align}
\end{lem}
{\it Proof.} In view of $e^{ (t^{1/2\alpha}- \tau^{1/2\alpha} - (t-\tau)^{1/2\alpha}) |\xi|_1} \in M_p$  and \eqref{N7.4},
\begin{align} \label{N7.6}
\|\Delta_j e^{t^{1/2\alpha} \Lambda} \mathscr{A}_{2\alpha} f\|_{p}
  & \lesssim \int_0^t e^{-c(t-\tau)2^{2\alpha j}}
 \|\Delta_j  e^{\tau^{1/2\alpha} \Lambda} f(\tau)\|_{p} d\tau.
\end{align}
Using Young's inequality, for $\gamma_1\le \gamma$, we have
\begin{align} \label{N7.7}
\|\Delta_j e^{t^{1/2\alpha} \Lambda} \mathscr{A}_{2\alpha} f\|_{L^\gamma_{t\in I}L^p_x}
    \lesssim   2^{- 2\alpha j (1+1/\gamma -1/\gamma_1) }
 \|\Delta_j  e^{\tau^{1/2\alpha} \Lambda} f(\tau)\|_{L^{\gamma_1}_{\tau \in I}L^p_x }.
\end{align}
Taking $l^q$ norms in both sides of \eqref{N7.7}, we have the results, as desired. $\hfill\Box$

\begin{thm}\label{GNSTh1a} {\rm (Analyticity for GNS: $p<\infty$)}
Let $n\ge 2$, $ 1/2 <\alpha<1$, $1< p < \infty$, $1\le q\le \infty$.  Assume that $u_0 \in   \dot{B}^{n/p+1-2\alpha}_{p, q}$, $\gamma_\pm = 2\alpha/(2\alpha-1\pm \varepsilon)$, $0<\varepsilon \ll 1$.  There exists a $T_{\max}= T_{\max}(u_0)>0$ such that \eqref{1.1} has a unique solution  $ u \in  \widetilde{L}^{\gamma_\pm}_{\rm loc}(0,T_{\max};  e^{{t}^{1/2\alpha}\Lambda}\dot{B}^{n/p \pm \varepsilon}_{p, q})   \cap \widetilde{L}^{\infty}([0,T_{\max}); e^{{t}^{1/2\alpha}\Lambda} \dot{B}^{n/p+1-2\alpha}_{p, q} )$.   If $T_{\max}<\infty$, then
$$
 \|u\|_{\widetilde{L}^{\gamma_+}(0,T_{\max};  e^{t^{1/2\alpha} \Lambda}\dot{B}^{n/p +\varepsilon}_{p, q}) \cap  \widetilde{L}^{\gamma_{-}}(0,T_{\max};  e^{t^{1/2\alpha} \Lambda}\dot{B}^{n/p-\varepsilon}_{p, q})}=\infty.
$$
 Moreover, if $u_0 \in   \dot{B}^{n/p+1-2\alpha}_{p, q}$ is sufficiently small, then $T_{\max}=\infty.$
\end{thm}

{\it Sketch of Proof.} Taking $\gamma= 2\alpha/(2\alpha-1)$, $\gamma_1 = \gamma/2$ in Lemmas \ref{7-1} and \ref{7-2}, we have for $I=[0,T)$,
\begin{align}
& \|U_{2\alpha}(t) u_0\|_{\widetilde{L}^{\gamma_\pm}(I; e^{t^{1/2\alpha}\Lambda} \dot{B}^{n/p\pm \varepsilon}_{p,
q})}\lesssim \|u_0 \|_{\dot{B}^{n/p+1-2\alpha}_{p, q}}, \label{N7.8} \\
& \|\nabla \mathscr{A}_{2\alpha} f\|_{\widetilde{L}^{\gamma_\pm}(I;
e^{t^{1/2\alpha} \Lambda} \dot{B}^{n/p\pm \varepsilon}_{p,q})}\lesssim
\|f\|_{\widetilde{L}^{\alpha/(2\alpha-1) }(I;e^{t^{1/2\alpha} \Lambda} \dot{B}^{n/p}_{p,q})}. \label{N7.9}
\end{align}
By \eqref{N7.9} and Lemma \ref{NSlem3},
\begin{align}
  \|\mathscr{A}_{2\alpha} \mathbb{P}\ \textrm{div}(u\otimes u) \|_{\widetilde{L}^{\gamma_\pm}(I;
e^{t^{1/2\alpha} \Lambda} \dot{B}^{n/p\pm \varepsilon}_{p,q})} & \lesssim
\| u\otimes u \|_{\widetilde{L}^{\alpha/(2\alpha-1) }(I;e^{t^{1/2\alpha} \Lambda} \dot{B}^{n/p}_{p,q})} \nonumber\\
& \lesssim
\| u \|^2_{\widetilde{L}^{\gamma_+}(I;  e^{t^{1/2\alpha} \Lambda}\dot{B}^{n/p +\varepsilon}_{p, q}) \cap  \widetilde{L}^{\gamma_{-}}(I;  e^{t^{1/2\alpha} \Lambda}\dot{B}^{n/p-\varepsilon}_{p, q}) }.  \label{N7.10}
\end{align}
Now, in view of \eqref{N7.8} and \eqref{N7.10}, we can use the same way as that of the proof of Theorem \ref{NSTh1a} to get the results, as desired. $\hfill\Box$

\begin{thm}\label{GNSTh1b} {\rm (Analyticity for GNS: $p=\infty$)}
Let $n\ge 2$, $ 1/2< \alpha < 1$, $1\le p \le \infty$.  Assume that $u_0 \in   M^{1-2\alpha}_{p, 1}$.  There exists a $T_{\max}= T_{\max}(u_0)>0$ such that \eqref{1.1} has a unique solution $u \in  \widetilde{L}^{2\alpha/(2\alpha-1)}_{\rm loc}(0,T_{\max}; {E}^{ct}_{p, 1})$ and $ (I-\Delta)^{(1-2\alpha)/2} u \in  \widetilde{L}^{\infty}([0,T_{\max}); {E}^{ct}_{p, 1})$. Moreover,
 if $T_{\max}<\infty$, then
 $$
   \|u\|_{\widetilde{L}^{2\alpha/(2\alpha-1)}(0,T_{\max}; {E}^{ct}_{p, 1})}=\infty.
 $$
\end{thm}

The proof of Theorem \ref{GNSTh1b} follows the ideas as in Theorems \ref{NSTh1a} and \ref{GNSTh1a} and we omit the details of the proof. \\

  \noindent{\bf Acknowledgment.} The authors are
supported in part by the National Science Foundation of China, grants 11201498 and  11271023, respectively.


\medskip
\footnotesize
\begin{thebibliography}{100}

\bibitem{BaBiTa12} H. Bae, A. Biswas and E. Tadmor,  Analyticity and decay estimates
of the Navier--Stokes equations in critical
Besov spaces, Arch. Rational Mech. Anal., {\bf 205} (2012), 963--991.

\bibitem{BenOk08} A. B\'{e}nyi, K.A. Okoudjou, Local well-posedness of nonlinear dispersive equations on modulation spaces, Bull. London Math. Soc., {\bf 41} (2009), 549--558.

\bibitem{BGOR} A. B\'{e}nyi,
K. Gr\"ochenig, K.A. Okoudjou and L.G. Rogers, Unimodular Fourier
multiplier for modulation spaces, J. Funct. Anal.,  \textbf{246} (2007), 366--384.

\bibitem{BL} J. Bergh and J. L\"{o}fstr\"{o}m,
 Interpolation Spaces,   Springer--Verlag,  1976.

\bibitem{BoPa08} J. Bourgain and N. Pavlovic,  Ill-posedness of the Navier--Stokes equations in a critical
space in 3D. J. Funct. Anal., $\mathbf{255}$ (2008), 2233--2247.

\bibitem{Can97} M. Cannone, A generalization of a theorem by Kato on Navier--Stokes equations, Rev. Mat. Iberoamericana, 13 (1997) 515--541.

\bibitem{CaPl96} M. Cannone and F. Planchon,  Self-similar solutions for Navier--Stokes equations in $\mathbf{R}^3$,  Comm. Partial Differ. Eq., $\mathbf{21}$ (1996), 179--193.

\bibitem{ChSh11} A. Cheskidov, and R. Shvydkoy, The regularity of weak solutions of the 3D Navier--Stokes equations in $B^{-1}_{\infty, \infty}$,  Arch. Rational Mech. Anal., {\bf 195} (2011), 159--169.

\bibitem{CH} J.-Y. Chemin,  { \rm Perfect incompressible fluids,  } Oxford University Press, Oxford,
(1998).

\bibitem{Ch99} J.-Y. Chemin,  Th\'{e}or\`{e}mes d'unicit\'{e} pour le syst\`{e}me de Navier--Stokes tridimensionnel,
J. Anal. Math., {\bf 77} (1999), 27--50.

\bibitem{ChFaSu12} J. Chen, D. Fan and L. Sun, Asymptotic estimates for unimodular Fourier multipliers on modulation spaces, Discrete Contin. Dyn. Syst.,  \textbf{32}  (2012),   467--485.

\bibitem{CorNik08} E. Cordero, F. Nicola,  Some new Strichartz estimates for the Schr\"odinger equation. J. Differential Equations, \textbf{245} (2008),   1945--1974.

 \bibitem{CorNik092} E. Cordero, F. Nicola,   Sharpness of some properties of Wiener amalgam and modulation spaces. Bull. Aust. Math. Soc., \textbf{80} (2009),  105--116.

\bibitem{DoLi09} H. Dong and D. Li, Optimal local smoothing and analyticity rate estimates for the generalized
Navier--Stokes equations. Commun. Math. Sci., {\bf 7} (2009), 67--80.

\bibitem{EsSeSv03} L. Escauriaza, G. Serigin and V. Sverak, $L_{3,\infty}$ solutions of Navier--Stokes equations
and backward uniquness, Uspekhi Mat. Nauk., \textbf{58} (2003), 3¨C44.

\bibitem{Fei83} H. G. Feichtinger, Modulation spaces on locally
compact Abelian group, Technical Report, University of Vienna, 1983.
Published in: ``Proc. Internat. Conf. on Wavelet and Applications",
99--140.  New Delhi Allied Publishers, India, 2003.

\bibitem{Fo97} C. Foias, What do the Navier--Stokes equations tell us about turbulence? Harmonic
Analysis and Nonlinear Differential Equations (Riverside, 1995). Contemp. Math., \textbf{208} (1997),
151--180.

\bibitem{FoTe89}  C. Foias and R. Temam,  Gevrey class regularity for the solutions of the Navier--Stokes
equations. J. Funct. Anal., \textbf{87} (1989), 359--369.

\bibitem{Ge08} P. Germain,  The second iterate for the Navier¨CStokes equation. J. Funct. Anal., $\mathbf{255}$ (2008), 2248--2264.

\bibitem{GePaSt07} P. Germain,  N. Pavlovic and G. Staffilani,  Regularity of solutions to the Navier¨C
Stokes equations evolving from small data in $BMO^{-1}$. Int. Math, Res. Not., $\mathbf{21}$ (2007), 35
pages.

\bibitem{Gr92} P. Gr\"obner, {\it Banachr\"aume Glatter Funktionen und
Zerlegungsmethoden}, Doctoral thesis, University of Vienna, 1992.

\bibitem{Groch01} K. Gr\"ochenig,  {\it Foundations of Time--Frequency
Analysis}, Birkh\"auser, Boston, MA, 2001.

\bibitem{GrKu98} Z. Gruji\v{c} and I. Kukavica, Space analyticity for the Navier--Stokes and related equations with initial data in $L^p$, J. Func. Anal., $\mathbf{152}$ (1998),
  247--466.

\bibitem{Ho60} L. H$\ddot {\rm o}$rmander,  { \rm Estimates for translation invariant
operators in $L^p$ spaces,  }   Acta Math.,  {\bf 104} (1960),
93--139.

\bibitem{Iw10} T. Iwabuchi, Navier--Stokes equations and nonlinear heat equations
in modulation spaces with negative derivative indices, J. Differential Equations, {\bf 248} (2010), 1972--2002.

\bibitem{Ka84} T. Kato, Strong $L^p$ solutions of the Navier--Stokes equations in $ \mathbb{{R}}^m$ with applications
to weak solutions,  Math. Z., \textbf{187} (1984), 471--480.

\bibitem{KoTa01} H. Koch, D. Tataru, Well-posedness for the Navier--Stokes equations, Adv. Math., {\rm 157} (2001) 22--35.

\bibitem{LeRi00}   P.G. Lemari\'{e}-Rieusset, {\rm On the analyticity of mild solutions
for the Navier--Stokes equations}, C. R. Acad. Sci. Paris,  Ser I,
{\bf 330} (2000), 183--186.

\bibitem{LeRi02} P.G. Lemari\'{e}-Rieusset,  Recent Developments in the Navier--Stokes Problem, Chapman
\& Hall/CRC Research Notes in Mathematics, vol \textbf{431},  Boca Raton, 2002.

\bibitem{Li} J.-L. Lions,  { \rm Quelques methods de resolution des probl$\grave{e}$mes
aux limites non lin$\acute{e}$aires,  }  (French)Paris:
Dunod/Gauthier-Villars, 1969.
\bibitem{LZ} P. Li and Z.  Zhai, Well-posedness and regularity of generalized Navier--Stokes equations in some critical Q--spaces, J. Funct. Anal.,  \textbf{259}  (2010), 2457--2519.

\bibitem{LY} P. Li and Q. Yang, Wavelets and the well-posedness of incompressible magneto-hydrodynamic equations in Besov type Q--space, J. Math. Anal. Appl.,  \textbf{405}  (2013),    661--686.

\bibitem{MiSa06} H. Miura and O. Sawada, On the regularizing rate estimates of Koch-Tataru's solution
to the Navier--Stokes equations,  Asymptot. Anal., \textbf{49} (2006), 1--15.

\bibitem{Pl96} F. Planchon, Global strong solutions in Sobolev or Lebesgue spaces to the incompressible Navier--Stokes equations in $\mathbb{R}^3$,
Ann. Inst. H. Poincare, AN, {\bf 13} (1996), 319--336.

\bibitem{SuTo07} M. Sugimoto amd N. Tomita, The dilation property of
modulation spaces and their inclusion relation with Besov spaces, J.
Funct. Anal., {\bf 248} (2007), 79--106.

\bibitem{To04} J. Toft, Continuity properties for modulation spaces,
with applications to pseudo-differential calculus, I,  J. Funct.
Anal., {\bf 207} (2004), 399--429.

\bibitem{Tr} H. Triebel,   Theory of Function Spaces,
  Birkh\"{a}user--Verlag,  1983.

  \bibitem{Wang02} B. Wang, {\rm The limit behavoir of solutions for the
complex Ginzburg--Landau  equation,} Commun. Pure Appl. Math.,
{\bf 55} (2002),  481--508.
\bibitem{Wang04} B. Wang, {\rm Exponential Besov spaces and their
applications to certain evolution equations with dissipation,}
Comm. Pure Appl. Anal., {\bf 3} (2004), 883--919.

\bibitem{WaHu07} B.   Wang and C.  Huang, Frequency-uniform decomposition method for the generalized BO, KdV and NLS equations,  J. Differential Equations, {\bf 239} (2007), 213--250.

\bibitem{WaHuHaGu11} B.  Wang, Z.   Huo, C.  Hao and Z.  Guo, {\it Harmonic Analysis Methods for Nonlinear Evolution Equations,} World Scientific, 2011.

\bibitem{Wa06}  B.  Wang, L.  Zhao, B. Guo,   {\rm  Isometric decomposition operators, function spaces $E^\lambda_{p,q}$ and their applications to nonlinear evolution equations,}  J. Funct. Anal., {\bf 233} (2006), 1--39.

\bibitem{WU3} J. Wu,  {\rm Lower bounds for an integral involving fractional Laplacians and the generalized
Navier-Stokes equations in Besov spaces},  Comm.Math.Phys., {\bf 263}
(2005),  803--831.

\bibitem{Yo10} T. Yoneda,  Ill-posedness of the 3D Navier--Stokes equations in a generalized Besov
space near $BMO^{-1}$. J. Funct. Anal., $\mathbf{258}$ (2010), 3376--3387.

\bibitem{YuZh12} X. Yu and Z. Zhai, Well-posedness for fractional Navier-Stokes equations in the largest critical spaces $\dot{B}^{-(2\beta-1)}_{\infty, \infty}(\mathbb{R}^n)$,  Math. Methods Appl. Sci.,  \textbf{35} (2012), 676--683.





\end{thebibliography}
\end{document}